\documentclass{article}
\usepackage{amsmath, amssymb, amsthm}
\usepackage{graphicx}
\usepackage{xcolor}
\usepackage{algorithm}
\usepackage[letterpaper, top=2.54cm, bottom=2.54cm, left=2.18cm, right=2.18cm, marginparwidth=1.75cm]{geometry}
\usepackage[colorlinks=true, linkcolor=blue, citecolor=blue, urlcolor=blue]{hyperref}
\newtheorem{theorem}{Theorem}
\newtheorem{lemma}{Lemma}
\newtheorem{assumption}{Assumption}
\newtheorem{remark}{Remark}
\newtheorem{definition}{Definition}

\title{Certified Nash Equilibrium Learning for Nonzero-Sum Jump-Diffusion Games: Policy Iteration and Off-Policy Q-Learning} 
\author{
    Pengwei Wang\footnote{School of Statistics and Mathematics, Shandong University of Finance and Economics, Jinan 250014, China. Email: \url{231113004@stu.sdufe.edu.cn}.},\quad
    Na Li\footnote{Corresponding author. School of Mathematical Sciences, Dalian University of Technology, Dalian 116024, China. Email: \url{lina2025@dlut.edu.cn}. The author acknowledges financial support from the NSFC (No. 12571475, No. 12171279).},\quad
    Jiacheng Fan\footnote{Department of Applied Mathematics, The Hong Kong Polytechnic University, Hong Kong SAR, China. Email: \url{jiacheng.fan@polyu.edu.hk}.},\quad
    Tielong Shen\footnote{School of Control Science and Engineering, Dalian University of Technology, Dalian 116024, China. Email: \url{tetusin@dlut.edu.cn}.},\quad
    Zuo Quan Xu\footnote{Department of Applied Mathematics, The Hong Kong Polytechnic University, Hong Kong SAR, China. Email: \url{maxu@polyu.edu.hk}. The author acknowledges financial support from the NSFC (Nos. 12571475 and 12571517), Hong Kong RGC (GRF 15203423 and 15204622), the PolyU-SDU Joint Research Center on Financial Mathematics, the CAS AMSS-PolyU Joint Laboratory of Applied Mathematics, the Research Centre for Quantitative Finance (1-CE03), and internal grants from The Hong Kong Polytechnic University.}
}
\date{\today}
\begin{document}
\maketitle
\begin{abstract}
    This paper establishes two model-free reinforcement learning methods for infinite-horizon \(N\)-player linear-quadratic nonzero-sum stochastic differential games whose drift, Brownian diffusion, and Poisson jump amplitudes depend on all players' controls. Control dependence in the diffusion and jump terms couples the players' simultaneous feedback updates. A verifiable certificate ensures that the partially model-free policy iteration map is contractive on an invariant closed ball centered at an admissible initial feedback profile. All feedback iterates are mean-square stabilizing, and the algorithm achieves local geometric convergence of feedback gains and value matrices to those associated with a closed-loop Nash equilibrium. The term certified refers to a verifiable convergence certificate that guarantees mean-square stability of all policy iterates and local geometric convergence to a closed-loop Nash equilibrium. To implement the policy updates from data, the off-policy Q-learning algorithm reuses the same behavior-policy data across iterations without requiring the system coefficients or jump intensity. Under quantitative persistence of excitation and a common initialization, exact population regression reproduces the policy-iteration sequence. The same local convergence guarantee therefore applies to off-policy Q-learning at the population level. Simulation results demonstrate the effectiveness of the two proposed algorithms.
\end{abstract}

\newpage
\section{Introduction}
\label{sec:introduction}
Stochastic systems involving strategic interactions typically evolve under continuous fluctuations together with abrupt random events. Examples include equipment failures in production lines, insurance claim arrivals, and sharp market price swings. In such settings, the decisions of individual agents may influence not only the expected evolution of the state but also its diffusion and jump amplitudes. When the agents pursue different objectives, these interactions lead naturally to nonzero-sum stochastic differential games with controlled jump-diffusion dynamics. Computing feedback Nash equilibria is challenging because uncertainty and strategic coupling enter the equilibrium equations simultaneously. The problem becomes more difficult when the system coefficients and jump characteristics are unknown and equilibrium policies must be learned from trajectory data. A further challenge arises when learning these policies through policy iteration, as the players’ feedback updates are coupled, and the resulting policy iteration procedure may fail to converge globally.

Linear quadratic (LQ) models are analytically tractable due to their linear dynamics and quadratic costs. In nonzero-sum differential games with this structure, Nash equilibrium~\cite{nash1951non} is the standard solution concept. Starr and Ho~\cite{starr1969nonzero} studied how to characterize Nash strategies in nonzero-sum differential games. Ba{\c{s}}ar and Olsder~\cite{bacsar1998dynamic} developed the general theory of dynamic noncooperative games, including feedback information structures and LQ formulations. In closed-loop LQ games, coupled Riccati equations provide equilibrium characterizations~\cite{freiling2002global,sun2019linear}. For stochastic LQ games with Poisson jumps, Zhang et al.~\cite{zhang2016poisson} characterized nonzero-sum Nash equilibria through coupled Riccati equations. Related developments include partially observed games with jumps~\cite{chen2014nonzero}, approximate Nash equilibria for controlled jump-diffusion mean-field interactions~\cite{benazzoli2019epsilon}, and stochastic Stackelberg LQ games with jumps~\cite{moon2021linear}. These results provide analytical foundations for equilibrium analysis, but their direct numerical implementation generally requires knowledge of the underlying dynamics.

Reinforcement learning (RL) and adaptive dynamic programming learn policies directly from observed state-control trajectories, avoiding pre-estimation of all system coefficients~\cite{sutton1998reinforcement,werbos1992approximate}. Watkins and Dayan~\cite{watkins1992q} introduced Q-learning, which learns an action-value function without first identifying the transition model. For continuous-time linear systems, Jiang and Jiang~\cite{jiang_computational_2012} developed data-driven policy iteration without requiring the system matrices. Subsequent work extended this approach to stochastic control with multiplicative noise~\cite{bian2016adaptive,li2022stochasticrl,zhang2024data} and to linear quadratic Gaussian (LQG) mean-field social control~\cite{xu2025social,xu2025robust}. More recent studies have begun to treat jumps explicitly. Guo et al.~\cite{guo2023reinforcement} studied finite-horizon learning for an unknown linear jump-diffusion system, and Gao et al.~\cite{gao2026reinforcement} developed policy-evaluation and Q-learning methods for entropy-regularized control of general jump diffusions. These results show that trajectory data can reduce model requirements, yet most works focus on single-agent or cooperative optimization instead of the coupled \(N\)-player Nash games.

Early RL methods for LQ nonzero-sum differential games mainly addressed deterministic dynamics. Vrabie and Lewis~\cite{vrabie2010integral} developed integral RL for two-player continuous-time games, and Vamvoudakis~\cite{vamvoudakis2015nashq} developed model-free Nash Q-learning for unknown deterministic continuous-time linear systems.
Li et al.~\cite{li2019policy} developed off-policy Q-learning for deterministic LQ nonzero-sum games using reusable behavior-policy data. Learning-based approaches have also been investigated for LQG mean-field games~\cite{xu2023modelfree,xu2025datadriven} and stochastic nonzero-sum games~\cite{zhang2024online,tian2024stochastic}. Recent advances in RL have also addressed games with jumps. Lu et al.~\cite{lu2025multiagent} studied multiagent investment games in jump-diffusion markets, characterized Nash equilibria under several utility functions, and developed parallel actor--critic learning. Hu et al.~\cite{hu2026entropy} formulated entropy-regularized learning for regime-switching jump-diffusion zero-sum games and proposed actor--critic policy improvement for investment-game settings. While these studies demonstrate that RL and jump-diffusion games can be combined, the infinite-horizon finite-player LQ nonzero-sum game with control-dependent diffusion and jumps requires further study.

There are two technical difficulties in learning feedback Nash equilibria for infinite-horizon stochastic LQ nonzero-sum games with Poisson jumps and unknown dynamics. First, when all players’ controls enter the diffusion and jump terms, each player’s stationarity equation involves the other players’ feedback gains. Simultaneous policy improvement therefore requires solving a coupled block system whose coefficients depend on the players’ value matrices. When the system coefficients and jump intensity are unknown, the data-driven method must recover both the individual value functions and the blocks needed for policy improvement from trajectory data. Second, convergence cannot be inferred directly from single-agent LQ control. Under standard assumptions, exact policy iteration in single-agent LQ control produces monotonically nonincreasing value matrices and converges to the optimal feedback from any stabilizing initial gain~\cite{zhang2024data}. In a nonzero-sum game, each player’s policy update can change the other players’ costs, so simultaneous updates can not decrease every player’s value matrix. The coupled algebraic Riccati equations may admit multiple stabilizing solutions or none~\cite{engwerda2005lq}, and the associated policy-iteration map is not globally contractive. Moreover, simultaneous and sequential updates define different maps, so convergence must be analyzed for the update rule used.

This paper addresses these difficulties through a partially model-free RL framework and its totally model-free off-policy Q-learning realization. We consider an infinite-horizon \(N\)-player LQ nonzero-sum stochastic differential game (LQ NZSSDG) driven jointly by a Brownian motion and a Poisson random measure, in which every player's control enters all three channels of the dynamics. A verifiable certificate is provided for local convergence of simultaneous policy iteration on a ball centered at a selected admissible initial feedback profile, without requiring prior knowledge of the equilibrium. The resulting policy iteration sequence serves as a reference for off-policy Q-learning, which reuses the same behavior-policy data across iterations without requiring the system coefficients or jump intensity. Under quantitative persistence of excitation, exact population regression reproduces this reference sequence, so the local convergence guarantee for policy iteration also applies to off-policy Q-learning at the population level.

The main contributions of this paper are summarized in the following three aspects.
\begin{enumerate}
	\item \textbf{Partially model-free policy-iteration RL with coupled diffusion and jump effects.} We develop a partially model-free policy-iteration RL algorithm that retains the cross-control terms induced by control-dependent diffusion and jump amplitudes. Integral Bellman policy evaluation is proved equivalent to Lyapunov recursions, and all players’ feedback gains are updated simultaneously through a coupled block system.
	
	\item \textbf{Local geometric convergence with an initial-policy-centered certificate.} We derive the Fréchet derivative of the policy-iteration map and establish local geometric convergence to the unique fixed point in an invariant ball, which we verify as a closed-loop Nash equilibrium. We also provide a verifiable certificate that checks mean-square stability, nonsingularity of the policy-improvement matrix, contraction, and a bound on the initial update residual.
	
	\item \textbf{Totally Model-free off-policy Q-learning with reusable trajectory data.} We derive an off-policy Q-learning algorithm with unknown system coefficients and jump intensity. At each iteration, the algorithm reuses the same batch of trajectory data collected under a fixed behavior policy with exploration noise. Under quantitative persistence of excitation, exact population regression uniquely recovers the Q-function matrices and reproduces the policy-iteration sequence, yielding the same local geometric convergence guarantee.
\end{enumerate}

The rest of this paper is organized as follows. Section~\ref{sec:game-formulation} characterizes feedback Nash equilibria. Section~\ref{sec:reinforcement-learning} develops a partially model-free RL algorithm and its off-policy Q-learning realization, and analyzes the convergence of both. Sections~\ref{sec:numerical-simulation} and~\ref{sec:conclusion} present the numerical example and conclusions.

\paragraph{Notation.}
\label{par:notation}
Denote by \(\mathbb{R}^{n}\) the \(n\)-dimensional Euclidean space and by \(\mathbb{R}^{n\times m}\) the set of real \(n\times m\) matrices. The Frobenius norm is \(\|A\|_F = (\operatorname{tr}(A^{\top} A))^{1/2}\), while \(\|v\|_2\) and \(\|M\|_2\) represent the Euclidean norm of vector \(v\) and the spectral norm of matrix \(M\), respectively. \(A^{\top}\) denotes transpose, \(\operatorname{vec}(\cdot)\) denotes columnwise vectorization, and \(\operatorname{sym}(M):= M+M^{\top}\). For \(x \in \mathbb{R}^{n}\), put \(\bar{x}:= [x_{1}^{2},x_{1}x_{2},\ldots,x_{1}x_{n},x_{2}^{2},\ldots,x_{n}^{2}]^{\top}\). Let \(\mathcal{S}^{n}\), \(\mathcal{S}_{+}^{n}\), and \(\mathcal{S}_{++}^{n}\) denote the symmetric, positive-semidefinite, and positive-definite matrices, respectively. All processes are defined on a complete filtered probability space \((\Omega,\mathcal{F},\mathbb{P},\mathbb{F})\) satisfying the usual conditions. The filtration \(\mathbb{F} = \{\mathcal{F}_{t}\}_{t \geq 0}\) is generated by two mutually independent stochastic processes and augmented by the collection of all \(\mathbb{P}\)-null sets. One is a standard one-dimensional Brownian motion \(\{W(t)\}_{t \geq 0}\), and the other is a Poisson random measure \(N\) defined on \(\mathbb{R}_{+}\times\mathcal{E}\), where \(\mathcal{E} = \mathbb{R}\setminus\{0\}\). The compensator of \(N\) is \(\overline{N}(\mathrm{d}t,\mathrm{d}e) = \lambda(\mathrm{d}e)\,\mathrm{d}t\), and the compensated Poisson random measure is \(\widetilde{N}(\mathrm{d}t,\mathrm{d}e):= N(\mathrm{d}t,\mathrm{d}e)-\lambda(\mathrm{d}e)\,\mathrm{d}t\). Hence, \(\{\widetilde{N}((0,t]\times S)\}_{t \geq 0}\) is a martingale for every \(S \in \mathcal{B}(\mathcal{E})\) with \(\lambda(S) < \infty\). Here, \(\lambda\) is a given \(\sigma\)-finite measure on \((\mathcal{E},\mathcal{B}(\mathcal{E}))\) such that \(\int_{\mathcal{E}}(1\wedge e^{2})\lambda(\mathrm{d}e) < \infty\). Moreover, \(L^{\lambda,2}(\mathcal{E};\mathbb{H})\) denotes the space of all \(\mathbb{H}\)-valued measurable functions \(r\) satisfying \(\int_{\mathcal{E}}\|r(e)\|_{\mathbb{H}}^{2}\lambda(\mathrm{d}e) < \infty\). Finally, \(L_{\mathbb{F}}^{2}([0,\infty);\mathbb{R}^{d})\) is the set of \(\mathbb{F}\)-progressively measurable processes \(\varphi\) satisfying \(\mathbb{E}\int_{0}^{\infty}\|\varphi(s)\|^{2}\mathrm{d}s < \infty\).

\section{Nash Equilibrium for \texorpdfstring{\(N\)}{N}-Player Nonzero-Sum Games on an Infinite Horizon}
\label{sec:game-formulation}
In this section, we formulate the infinite-horizon jump-diffusion game and present conditions for a closed-loop Nash equilibrium. Consider the following continuous-time stochastic differential equation with Poisson jumps:
\begin{equation}
    \label{eq:game-dynamics}
    \left\{
    \begin{aligned}
        \mathrm{d}X(s) & =\left[AX(s)+\sum_{j = 1}^{N}B_{j}u_{j}(s)\right]\mathrm{d}s
        +\left[CX(s)+\sum_{j = 1}^{N}D_{j}u_{j}(s)\right]\mathrm{d}W(s)\\
        & \quad+\int_{\mathcal{E}}\left[E(e)X(s-)+\sum_{j = 1}^{N}F_{j}(e)u_{j}(s)\right]\widetilde{N}(\mathrm{d}s,\mathrm{d}e),\quad s \geq 0,\\
        X(0) & =x \in \mathbb{R}^n,
    \end{aligned}
    \right.
\end{equation}
where \(X(s) \in \mathbb{R}^{n}\) is the state process, \(u_{j}(s) \in \mathbb{R}^{m_{j}}\) is the control input of player \(j \in \mathbb N\), with \(\mathbb N = \{1,\ldots,N\}\). The matrices \(A,C \in \mathbb{R}^{n\times n}\), \(B_{j},D_{j} \in \mathbb{R}^{n\times m_{j}}, j \in \mathbb N\) are constant system matrices, while \(E(\cdot) \in L^{\lambda,2}(\mathcal{E};\mathbb{R}^{n\times n}), F_{j}(\cdot) \in L^{\lambda,2}(\mathcal{E};\mathbb{R}^{n\times m_{j}}), j \in \mathbb N\), are deterministic matrix-valued functions.

\begin{definition}
    \label{def:mean-square-stabilizability}
    System~\eqref{eq:game-dynamics} is called mean-square stabilizable if there exist \(N\) constant matrices \(K_j\in\mathbb{R}^{m_j\times n}\), \(j\in\mathbb N\), such that, for every initial state \(x\in\mathbb{R}^n\), the solution of
    \begin{equation*}
        \label{eq:game-closed-loop-dynamics}
        \left\{
        \begin{aligned}
            \mathrm{d}X(s)
            & =\left(A+\sum_{j = 1}^{N}B_{j}K_{j}\right)X(s)\,\mathrm{d}s
            +\left(C+\sum_{j = 1}^{N}D_{j}K_{j}\right)X(s)\,\mathrm{d}W(s)\\
            & \quad+\int_{\mathcal{E}}\left(E(e)+\sum_{j = 1}^{N}F_{j}(e)K_{j}\right)X(s-)\,\widetilde{N}(\mathrm{d}s,\mathrm{d}e),\quad s \geq 0,\\
            X(0) & =x,
        \end{aligned}
        \right.
    \end{equation*}
    satisfies \(\lim_{s \rightarrow \infty}\mathbb{E}\bigl[X(s)^{\top}X(s)\bigr] = 0\). In this case, \((K_{1},\ldots,K_{N})\) is called a stabilizer of the system and \((K_{1}X(\cdot),\ldots,K_{N}X(\cdot))\) is called stabilizing for the initial state \(x\).
\end{definition}

\begin{assumption}
    \label{ass:mean-square-stabilizability}
    The system~\eqref{eq:game-dynamics} is mean-square stabilizable.
\end{assumption}

\begin{lemma}
    \label{lem:mean-square-stability}
    Let \(K_{j} \in \mathbb{R}^{m_{j}\times n}\), \(j \in \mathbb N\). Then \((K_{1},\ldots,K_{N})\) is a mean-square stabilizer of~\eqref{eq:game-dynamics} if and only if there exists \(P \in \mathcal{S}_{++}^{n}\) such that
    \begin{equation*}
        \label{eq:game-stability-inequality}
        \begin{aligned}
            & \left(A+\sum_{j = 1}^{N}B_{j}K_{j}\right)^{\top}P
            +P\left(A+\sum_{j = 1}^{N}B_{j}K_{j}\right)
            +\left(C+\sum_{j = 1}^{N}D_{j}K_{j}\right)^{\top}
            P\left(C+\sum_{j = 1}^{N}D_{j}K_{j}\right) \\
            & \quad
            +\int_{\mathcal{E}}
            \left(E(e)+\sum_{j = 1}^{N}F_{j}(e)K_{j}\right)^{\top}
            P
            \left(E(e)+\sum_{j = 1}^{N}F_{j}(e)K_{j}\right)
            \lambda(\mathrm{d}e)
            < 0.
        \end{aligned}
    \end{equation*}
    In this case, for any \(\Lambda \in \mathcal{S}^{n}\) (respectively, \(\mathcal{S}_{+}^{n}\) and \(\mathcal{S}_{++}^{n}\)), the Lyapunov equation
    \begin{equation*}
        \label{eq:game-generalized-lyapunov}
        \begin{aligned}
            & \left(A+\sum_{j = 1}^{N}B_{j}K_{j}\right)^{\top}P
            +P\left(A+\sum_{j = 1}^{N}B_{j}K_{j}\right)
            +\left(C+\sum_{j = 1}^{N}D_{j}K_{j}\right)^{\top}
            P\left(C+\sum_{j = 1}^{N}D_{j}K_{j}\right) \\
            & \quad
            +\int_{\mathcal{E}}
            \left(E(e)+\sum_{j = 1}^{N}F_{j}(e)K_{j}\right)^{\top}
            P
            \left(E(e)+\sum_{j = 1}^{N}F_{j}(e)K_{j}\right)
            \lambda(\mathrm{d}e)
            +\Lambda = 0
        \end{aligned}
    \end{equation*}
    admits a unique solution \(P \in \mathcal{S}^{n}\) (respectively, \(\mathcal{S}_{+}^{n}\) and \(\mathcal{S}_{++}^{n}\)).
\end{lemma}

If the system is mean-square stabilizable for a given \(x \in \mathbb{R}^{n}\), define the corresponding set of admissible controls as
\begin{equation*}
    \label{eq:game-admissible-controls}
    \begin{aligned}
        \mathcal{U}_{ad}^{x}:= \Biggl\{(u_{1}(\cdot),\ldots,u_{N}(\cdot)):
        & \ u_{j}(\cdot) \in L_{\mathbb{F}}^{2}([0,\infty);\mathbb{R}^{m_{j}})\ \text{with}\ j \in \mathbb N\\
        & \ \text{and}\ (u_{1}(\cdot),u_{2}(\cdot),\ldots,u_{N}(\cdot))\ \text{is stabilizing w.r.t}\ x \Biggr\}.
    \end{aligned}
\end{equation*}
For simplicity, denote \(\mathbf{u}(\cdot) = (u_{1}(\cdot),\ldots,u_{N}(\cdot))\). Define the cost functional for player \(i\), \(i = 1,2,\ldots,N\), as
\begin{equation}
    \label{eq:game-cost-functional}
    J_{i}(t,x;\mathbf{u}(\cdot))
    = \mathbb{E}^{\mathcal{F}_{t}}\int_{t}^{\infty}
    \left[X(s)^{\top} N_{i}X(s)+\sum_{j = 1}^{N}u_{j}(s)^{\top} R_{ij}u_{j}(s)\right]\mathrm{d}s,
\end{equation}
where \(N_{i} \in \mathcal{S}_{++}^{n}\) and \(R_{ij} \in \mathcal{S}_{++}^{m_{j}}\) for all \(i,j \in \mathbb N\). When the initial time is zero, we write \(J_{i}(x;\mathbf{u}(\cdot)):= J_{i}(0,x;\mathbf{u}(\cdot))\).

\textbf{Problem (LQ NZSSDG).} For any given \(x \in \mathbb{R}^{n}\), the problem is to find an admissible control profile \(\mathbf{u}^{*}(\cdot) \in \mathcal{U}_{ad}^{x}\) such that, for each \(i \in \mathbb N\),
\begin{equation*}
    \label{eq:game-nash-inequality}
    J_{i}(x;\mathbf{u}^{*}(\cdot)) \leq J_{i}(x;u_{i}(\cdot), u_{-i}^{*}(\cdot)), \qquad
    \forall\,u_{i}(\cdot)\ \text{s.t.}\ (u_{i}(\cdot),u_{-i}^{*}(\cdot)) \in \mathcal{U}_{ad}^{x},
\end{equation*}
where \(u_{-i}^{*} = (u_{1}^{*},\dots, u_{i-1}^{*}, u_{i+1}^{*}, \dots, u_{N}^{*})\). If \(\mathbf{u}^{*}(\cdot)\) further admits a linear state feedback representation
\begin{equation*}
    \label{eq:game-nash-feedback}
    \mathbf{u}^{*}(\cdot) = (K_{1}^{*}X(\cdot),\ldots,K_{N}^{*}X(\cdot)) \in \mathcal{U}_{ad}^{x},
\end{equation*}
then it is called a \emph{feedback Nash equilibrium} of Problem (LQ NZSSDG).

Define the Hamiltonian function associated with~\eqref{eq:game-dynamics} and~\eqref{eq:game-cost-functional} as
\begin{equation*}
    \label{eq:game-hamiltonian}
    \begin{aligned}
        \mathcal{H}_{i}(x,u_{1},\ldots,u_{N};P_{i})
        & =
        x^{\top} N_{i}x+\sum_{j = 1}^{N}u_{j}^{\top} R_{ij}u_{j}\\
        & \quad
        +x^{\top} P_{i}\left(Ax+\sum_{j = 1}^{N}B_{j}u_{j}\right)
        +\left(Ax+\sum_{j = 1}^{N}B_{j}u_{j}\right)^{\top} P_{i}x\\
        & \quad
        +\left(Cx+\sum_{j = 1}^{N}D_{j}u_{j}\right)^{\top}
        P_{i}
        \left(Cx+\sum_{j = 1}^{N}D_{j}u_{j}\right)\\
        & \quad+\int_{\mathcal{E}}\left(E(e)x+\sum_{j = 1}^{N}F_{j}(e)u_{j}\right)^{\top}
        P_{i}\left(E(e)x+\sum_{j = 1}^{N}F_{j}(e)u_{j}\right)\lambda(\mathrm{d}e).
    \end{aligned}
\end{equation*}

Before designing iterative algorithms, we need a tractable characterization of feedback Nash equilibria. In the LQ setting, the value function of each player is quadratic, and the equilibrium conditions reduce to a coupled system consisting of a Lyapunov equation for each player's value matrix and a stationarity condition for each player's gain.

\begin{lemma}[Equilibrium characterization]
    \label{lem:nash-characterization}
    Let \(\mathbf{u}^{*}(\cdot) = (K_{1}^{*}X(\cdot),\ldots,K_{N}^{*}X(\cdot)) \in \mathcal{U}_{ad}^{x}\) be a linear state-feedback strategy. Then \(\mathbf{u}^{*}(\cdot)\) is a closed-loop Nash equilibrium of Problem (LQ NZSSDG) if and only if, for each \(i = 1,\ldots,N\), the solution \(P_{i} \in \mathcal{S}_{++}^{n}\) of the Lyapunov equations
    \begin{equation}
        \label{eq:game-nash-lyapunov}
        \begin{aligned}
            & \left(A+\sum_{j = 1}^{N}B_{j}K_{j}^{*}\right)^{\top} P_{i}
            +P_{i}\left(A+\sum_{j = 1}^{N}B_{j}K_{j}^{*}\right)
            +\left(C+\sum_{j = 1}^{N}D_{j}K_{j}^{*}\right)^{\top}
            P_{i}
            \left(C+\sum_{j = 1}^{N}D_{j}K_{j}^{*}\right)\\
            & \quad+\int_{\mathcal{E}}\left(E(e)+\sum_{j = 1}^{N}F_{j}(e)K_{j}^{*}\right)^{\top}
            P_{i}\left(E(e)+\sum_{j = 1}^{N}F_{j}(e)K_{j}^{*}\right)\lambda(\mathrm{d}e)
            +N_{i}+\sum_{j = 1}^{N}(K_{j}^{*})^{\top} R_{ij}K_{j}^{*}
            = 0
        \end{aligned}
    \end{equation}
    satisfies the following conditions:
    \begin{equation}
        \label{eq:game-nash-stationarity}
        \begin{aligned}
            R_{ii} & +D_{i}^{\top} P_{i}D_{i}+\int_{\mathcal{E}}F_{i}(e)^{\top} P_{i}F_{i}(e)\lambda(\mathrm{d}e) > 0,\\
            B_{i}^{\top} P_{i}+D_{i}^{\top} P_{i} & \left(C+\sum_{j = 1}^{N}D_{j}K_{j}^{*}\right)
            +\int_{\mathcal{E}}F_{i}(e)^{\top} P_{i}\left(E(e)+\sum_{j = 1}^{N}F_{j}(e)K_{j}^{*}\right)\lambda(\mathrm{d}e)
            +R_{ii}K_{i}^{*} = 0.
        \end{aligned}
    \end{equation}
\end{lemma}

\begin{proof}
    \textbf{Sufficiency.} Fix an arbitrary player \(i \in \{1,\ldots,N\}\) and fix the other players' equilibrium strategies \(u_{j}^{*}(s) = K_{j}^{*}X(s)\) for all \(j \neq i\). Let player \(i\) use an arbitrary admissible deviation \(u_{i}(\cdot)\). Write it as
    \begin{equation}
        \label{eq:game-control-deviation}
        u_{i}(s) = K_{i}^{*}X(s)+v(s),
    \end{equation}
    where \(v(\cdot)\) is an adapted process such that the resulting cost is finite. Under~\eqref{eq:game-control-deviation}, the state system~\eqref{eq:game-dynamics} becomes
    \begin{equation*}
        \label{eq:game-deviation-dynamics}
        \begin{aligned}
            \mathrm{d}X(s) & =\bigl(A_{K}^{*}X(s)+B_{i}v(s)\bigr)\mathrm{d}s
            +\bigl(C_{K}^{*}X(s)+D_{i}v(s)\bigr)\mathrm{d}W(s)\\
            & \quad+\int_{\mathcal{E}}\bigl(E_{K}^{*}(e)X(s-)+F_{i}(e)v(s)\bigr)\widetilde{N}(\mathrm{d}s,\mathrm{d}e),
        \end{aligned}
    \end{equation*}
    where
    \begin{equation*}
        \label{eq:game-equilibrium-coefficients}
        A_{K}^{*}:= A+\sum_{j = 1}^{N} B_{j}K_{j}^{*}, \qquad
        C_{K}^{*}:= C+\sum_{j = 1}^{N} D_{j}K_{j}^{*},\qquad
        E_{K}^{*}(e):= E(e)+\sum_{j = 1}^{N}F_{j}(e)K_{j}^{*}.
    \end{equation*}
    Applying the It\^o formula with jumps to \(X(s)^{\top} P_{i}X(s)\), we obtain
    \begin{equation}
        \label{eq:game-deviation-ito}
        \begin{aligned}
            \mathrm{d}(X^{\top} P_{i}X)
            & =
            \Bigl[
                X^{\top}\Bigl((A_{K}^{*})^{\top} P_{i}+P_{i}A_{K}^{*}+(C_{K}^{*})^{\top} P_{i}C_{K}^{*}
                +\int_{\mathcal{E}}(E_{K}^{*}(e))^{\top} P_{i}E_{K}^{*}(e)\lambda(\mathrm{d}e)\Bigr)X\\
                & \qquad
                +2v^{\top}\Bigl(B_{i}^{\top} P_{i}+D_{i}^{\top} P_{i}C_{K}^{*}
                +\int_{\mathcal{E}}F_{i}(e)^{\top} P_{i}E_{K}^{*}(e)\lambda(\mathrm{d}e)\Bigr)X\\
                & \qquad+v^{\top}\Bigl(D_{i}^{\top} P_{i}D_{i}
                +\int_{\mathcal{E}}F_{i}(e)^{\top} P_{i}F_{i}(e)\lambda(\mathrm{d}e)\Bigr)v
                \Bigr]\mathrm{d}s\\
            & \quad
            +\{\cdots\}\,\mathrm{d}W(s) + \{\cdots\}\,\widetilde{N}(\mathrm{d}s,\mathrm{d}e).
        \end{aligned}
    \end{equation}
    Integrating~\eqref{eq:game-deviation-ito} over \([0,T]\) and taking expectations yields
    \begin{equation}
        \label{eq:game-deviation-expectation}
        \begin{aligned}
            \mathbb{E}\bigl[X(T)^{\top} P_{i}X(T)\bigr]-x^{\top} P_{i}x
            & =\mathbb{E}\int_{0}^{T} \Biggl[X^{\top}\Bigl((A_{K}^{*})^{\top} P_{i}+P_{i}A_{K}^{*}\\
                & \qquad +(C_{K}^{*})^{\top} P_{i}C_{K}^{*}
                +\int_{\mathcal{E}}(E_{K}^{*}(e))^{\top} P_{i}E_{K}^{*}(e)\,\lambda(\mathrm{d}e)\Bigr)X \\
                & \qquad +2v^{\top}\left(B_{i}^{\top} P_{i}+D_{i}^{\top} P_{i}C_{K}^{*}
                +\int_{\mathcal{E}}F_{i}^{\top} P_{i}E_{K}^{*}(e)\,\lambda(\mathrm{d}e)\right)X \\
                & \qquad +v^{\top}\left(D_{i}^{\top} P_{i}D_{i}
                +\int_{\mathcal{E}}F_{i}^{\top} P_{i}F_{i}\,\lambda(\mathrm{d}e)\right)v\Biggr]\mathrm{d}s,
        \end{aligned}
    \end{equation}
    Substitute the Lyapunov equation~\eqref{eq:game-nash-lyapunov} into~\eqref{eq:game-deviation-expectation}:
    \begin{equation}
        \label{eq:game-deviation-lyapunov-identity}
        \begin{aligned}
            x^{\top} P_{i}x-\mathbb{E}\bigl[X(T)^{\top} P_{i}X(T)\bigr]
            & =\mathbb{E}\int_{0}^{T} \Biggl[X^{\top} N_{i}X+\sum_{j = 1}^{N} X^{\top} (K_{j}^{*})^{\top} R_{ij}K_{j}^{*}X \\
                & \qquad -2v^{\top}\left(B_{i}^{\top} P_{i}+D_{i}^{\top} P_{i}C_{K}^{*}
                +\int_{\mathcal{E}}F_{i}^{\top} P_{i}E_{K}^{*}(e)\,\lambda(\mathrm{d}e)\right)X\\
                & \qquad-v^{\top}\left(D_{i}^{\top} P_{i}D_{i}
                +\int_{\mathcal{E}}F_{i}^{\top} P_{i}F_{i}\,\lambda(\mathrm{d}e)\right)v\Biggr]\mathrm{d}s.
        \end{aligned}
    \end{equation}
    Under the deviation~\eqref{eq:game-control-deviation} and \(u_{j}^{*}(s) = K_{j}^{*}X(s)\) for \(j \neq i\), the running cost satisfies
    \begin{equation*}
        \label{eq:game-deviation-running-cost}
        \begin{aligned}
            X^{\top} N_{i}X+\sum_{j = 1}^{N} u_{j}^{\top} R_{ij}u_{j} = X^{\top} N_{i}X+\sum_{j = 1}^{N} X^{\top} (K_{j}^{*})^{\top} R_{ij}K_{j}^{*}X + 2v^{\top} R_{ii}K_{i}^{*}X+v^{\top} R_{ii}v.
        \end{aligned}
    \end{equation*}
    Combining this identity with~\eqref{eq:game-deviation-lyapunov-identity}, we obtain
    \begin{equation}
        \label{eq:game-finite-horizon-cost}
        \begin{aligned}
            \mathbb{E}\int_{0}^{T}
            & \left(X^{\top} N_{i}X+\sum_{j = 1}^{N} u_{j}^{\top} R_{ij}u_{j}\right)\mathrm{d}s = x^{\top} P_{i}x-\mathbb{E}\bigl[X(T)^{\top} P_{i}X(T)\bigr] \\
            & \quad +\mathbb{E}\int_{0}^{T} \Bigl[
                2v^{\top}\left(R_{ii}K_{i}^{*}+B_{i}^{\top} P_{i}+D_{i}^{\top} P_{i}C_{K}^{*}
                +\int_{\mathcal{E}}F_{i}^{\top} P_{i}E_{K}^{*}(e)\,\lambda(\mathrm{d}e)\right)X\\
                & \qquad+v^{\top}\left(R_{ii}+D_{i}^{\top} P_{i}D_{i}
                +\int_{\mathcal{E}}F_{i}^{\top} P_{i}F_{i}\,\lambda(\mathrm{d}e)\right)v
                \Bigr]\mathrm{d}s.
        \end{aligned}
    \end{equation}
    By the stationarity condition in~\eqref{eq:game-nash-stationarity},~\eqref{eq:game-finite-horizon-cost} reduces to
    \begin{equation*}
        \label{eq:game-finite-horizon-deviation}
        \begin{aligned}
            \mathbb{E}\int_{0}^{T}
            \left(X^{\top} N_{i}X+\sum_{j = 1}^{N} u_{j}^{\top} R_{ij}u_{j}\right)\mathrm{d}s
            & =x^{\top} P_{i}x-\mathbb{E}\bigl[X(T)^{\top} P_{i}X(T)\bigr]\\
            & \quad +\mathbb{E}\int_{0}^{T} v^{\top} \left(R_{ii}+D_{i}^{\top} P_{i}D_{i}
            +\int_{\mathcal{E}}F_{i}^{\top} P_{i}F_{i}\,\lambda(\mathrm{d}e)\right) v \mathrm{d}s.
        \end{aligned}
    \end{equation*}
    Letting \(T \to \infty\), by admissibility and mean-square stability, we obtain
    \begin{equation*}
        \label{eq:game-infinite-horizon-deviation}
        \begin{aligned}
            J_{i}\bigl(x;u_{i},u_{-i}^{*}\bigr)
            & =
            x^{\top} P_{i}x
            +\mathbb{E}\int_{0}^{\infty}
            v^{\top}\left(R_{ii}+D_{i}^{\top} P_{i}D_{i}
            +\int_{\mathcal{E}}F_{i}^{\top} P_{i}F_{i}\,\lambda(\mathrm{d}e)\right)v\,\mathrm{d}s \geq x^{\top} P_{i}x
            = J_{i}(x;\mathbf{u}^{*}).
        \end{aligned}
    \end{equation*}
    Therefore, player \(i\) cannot reduce its cost by any unilateral admissible deviation when the other players use \(\{K_{j}^{*}\}_{j \neq i}\). Since \(i\) is arbitrary, \(\mathbf{u}^{*}\) is a closed-loop Nash equilibrium.

    \medskip
    \textbf{Necessity.} Conversely, suppose that \(\mathbf{u}^{*}(\cdot) = (K_{1}^{*}X(\cdot),\ldots,K_{N}^{*}X(\cdot))\) is a closed-loop Nash equilibrium of Problem (LQ NZSSDG). Fix an arbitrary player \(i \in \mathbb N\) and keep the strategies of all other players at equilibrium, i.e., \(u_{j}(s) = K_{j}^{*}X(s)\) for all \(j \neq i\). Then player \(i\) faces a stochastic linear quadratic control problem with state dynamics
    \begin{equation*}
        \label{eq:game-best-response-dynamics}
        \begin{aligned}
            \mathrm{d}X(s)
            & =\left[\left(A+\sum_{j \neq i}B_{j}K_{j}^{*}\right)X(s)+B_{i}u_{i}(s)\right]\mathrm{d}s\\
            & \quad+\left[\left(C+\sum_{j \neq i}D_{j}K_{j}^{*}\right)X(s)+D_{i}u_{i}(s)\right]\mathrm{d}W(s)\\
            & \quad+\int_{\mathcal{E}}\left[\left(E(e)+\sum_{j \neq i}F_{j}(e)K_{j}^{*}\right)X(s)+F_{i}(e)u_{i}(s)\right]\widetilde{N}(\mathrm{d}s,\mathrm{d}e).
        \end{aligned}
    \end{equation*}
    and cost functional
    \begin{equation*}
        \label{eq:game-best-response-cost}
        J_{i}(x;u_{i},K_{-i}^{*})
        = \mathbb{E}\!\left[\int_{0}^{\infty} \Bigl(X(s)^{\top}\!\bigl(N_{i}+\sum_{j \neq i}(K_{j}^{*})^{\top} R_{ij}K_{j}^{*}\bigr)X(s)+u_{i}(s)^{\top} R_{ii}u_{i}(s)\Bigr)\mathrm{d}s\right].
    \end{equation*}

    Since \(\mathbf{u}^{*}(\cdot)\) is a closed-loop Nash equilibrium, \(u_{i}^{*}(s) = K_{i}^{*}X(s)\) is optimal for this LQ problem. Its quadratic value function is written as \(V_{i}(x) = x^{\top} P_{i}x\). Based on the definition of the Hamiltonian function, the Hamiltonian function of player \(i\) can be written as \(\mathcal{H}_{i}(x,u_{i},u_{-i}^{*};P_{i})\). The value function satisfies the HJB equation
    \begin{equation*}
        \label{eq:game-hjb-equation}
        \min_{u_{i} \in \mathbb{R}^{m_{i}}}\mathcal{H}_{i}(x,u_{i},u_{-i}^{*};P_{i}) = 0.
    \end{equation*}
    Differentiating \(\mathcal{H}_{i}(x,u_{i},u_{-i}^{*};P_{i})\) term by term with respect to \(u_{i}\) yields
    \begin{equation}
        \label{eq:game-hamiltonian-gradient}
        \begin{aligned}
            \frac{\partial \mathcal{H}_{i}}{\partial u_{i}}
            & =2R_{ii}u_{i}+2B_{i}^{\top} P_{i}x\\
            & \quad+2 D_{i}^{\top} P_{i}\left((C+\sum_{j \neq i}D_{j}K_{j}^{*})x+D_{i}u_{i}\right)\\
            & \quad+2 \int_{\mathcal{E}}F_{i}(e)^{\top} P_{i}\left((E(e)+\sum_{j \neq i}F_{j}(e)K_{j}^{*})x+F_{i}(e)u_{i}\right)\lambda(\mathrm{d}e)\\
            & =2\left(B_{i}^{\top} P_{i}+D_{i}^{\top} P_{i}(C+\sum_{j \neq i}D_{j}K_{j}^{*})
            +\int_{\mathcal{E}}F_{i}(e)^{\top} P_{i}(E(e)+\sum_{j \neq i}F_{j}(e)K_{j}^{*})\lambda(\mathrm{d}e)\right)x \\
            & \quad+2\left(R_{ii}+D_{i}^{\top} P_{i}D_{i}
            +\int_{\mathcal{E}}F_{i}(e)^{\top} P_{i}F_{i}(e) \lambda(\mathrm{d}e)\right)u_{i}.
        \end{aligned}
    \end{equation}
    Since \(u_{i}^{*}(x) = K_{i}^{*}x\) is optimal, it minimizes \(\mathcal{H}_{i}\) with respect to \(u_{i}\). Thus, the first-order optimality condition holds:
    \begin{equation*}
        \label{eq:game-first-order-condition}
        \begin{aligned}
            \frac{\partial \mathcal{H}_{i}}{\partial u_{i}} & \Big|_{u_{i} = K_{i}^{*}x}\\
            & =2\left(B_{i}^{\top} P_{i}+D_{i}^{\top} P_{i} \left(C+\sum_{j = 1}^{N}D_{j}K_{j}^{*}\right) +\int_{\mathcal{E}}F_{i}(e)^{\top} P_{i}\left(E(e)+\sum_{j = 1}^{N}F_{j}(e)K_{j}^{*}\right)\lambda(\mathrm{d}e) +R_{ii}K_{i}^{*}\right)x = 0.
        \end{aligned}
    \end{equation*}
    Because this holds for all \(x\), the matrix in parentheses is zero, which is the second relation in~\eqref{eq:game-nash-stationarity}. Differentiating~\eqref{eq:game-hamiltonian-gradient} once more gives the Hessian
    \begin{equation*}
        \label{eq:game-hamiltonian-hessian}
        \frac{\partial^{2} \mathcal{H}_{i}}{\partial u_{i}^{2}}
        = 2\left(R_{ii}+D_{i}^{\top} P_{i}D_{i}+\int_{\mathcal{E}}F_{i}(e)^{\top} P_{i}F_{i}(e) \lambda(\mathrm{d}e)\right).
    \end{equation*}
    For \(u_{i}^{*}\) to be a unique minimizer, the Hamiltonian must be strictly convex in \(u_{i}\), which requires
    \begin{equation*}
        \label{eq:game-hamiltonian-convexity}
        R_{ii}+D_{i}^{\top} P_{i}D_{i}+\int_{\mathcal{E}}F_{i}(e)^{\top} P_{i}F_{i}(e) \lambda(\mathrm{d}e) > 0.
    \end{equation*}
    Substituting the optimal feedback \(u_{i}^{*}(x) = K_{i}^{*}x\) into the HJB equation and observing the coefficients of \(x^{\top}(\cdot)x\) on both sides, we obtain
    \begin{equation*}
        \label{eq:game-hjb-lyapunov}
        \begin{aligned}
            & \left(A+\sum_{j = 1}^{N}B_{j}K_{j}^{*}\right)^{\top} P_{i}
            +P_{i}\left(A+\sum_{j = 1}^{N}B_{j}K_{j}^{*}\right)
            +\left(C+\sum_{j = 1}^{N}D_{j}K_{j}^{*}\right)^{\top}
            P_{i}
            \left(C+\sum_{j = 1}^{N}D_{j}K_{j}^{*}\right)\\
            & \quad+\int_{\mathcal{E}}\left(E(e)+\sum_{j = 1}^{N}F_{j}(e)K_{j}^{*}\right)^{\top}
            P_{i}\left(E(e)+\sum_{j = 1}^{N}F_{j}(e)K_{j}^{*}\right)\lambda(\mathrm{d}e)
            +N_{i}+\sum_{j = 1}^{N}(K_{j}^{*})^{\top} R_{ij}K_{j}^{*}
            = 0.
        \end{aligned}
    \end{equation*}
    which is exactly the Lyapunov equation~\eqref{eq:game-nash-lyapunov}.
\end{proof}

\section{Reinforcement Learning Algorithms for \texorpdfstring{\(N\)}{N}-Player Nash Equilibrium}
\label{sec:reinforcement-learning}

In this section, we an develop online partially model-free RL algorithm and off-policy Q-learning for Problem (LQ NZSSDG). We establish local convergence and a sufficient certificate for simultaneous policy iteration, and prove population equivalence of the two algorithms under a quantitative excitation condition.

\subsection{Online partially model-free policy-iteration RL algorithm}
\label{subsec:policy-iteration}
In this subsection, we use the equilibrium characterization in Section~\ref{sec:game-formulation} to define policy evaluation and simultaneous policy improvement. We formulate an online partially model-free RL algorithm, establish local geometric convergence, and derive sufficient tests for its assumptions. The improvement step uses partial system information without the drift matrix \(A\).

Let \(\mathcal{K}:= \prod_{j = 1}^{N}\mathbb{R}^{m_{j}\times n}\). For \(K = (K_{1},\ldots,K_{N}) \in \mathcal{K}\), define
\begin{equation*}
    \label{eq:pi-closed-loop-coefficients}
    A_{K}:= A+\sum_{j = 1}^{N}B_{j}K_{j},\qquad
    C_{K}:= C+\sum_{j = 1}^{N}D_{j}K_{j},\qquad
    E_{K}(e):= E(e)+\sum_{j = 1}^{N}F_{j}(e)K_{j}.
\end{equation*}
Define the generalized Lyapunov operator
\begin{equation*}
    \label{eq:pi-lyapunov-operator}
    \mathcal{L}_{K}(P):=
    A_{K}^{\top} P+PA_{K}+C_{K}^{\top} PC_{K}+
    \int_{\mathcal{E}}E_{K}(e)^{\top} PE_{K}(e)\lambda(\mathrm{d}e),
\end{equation*}
and the running-cost matrix
\begin{equation*}
    \label{eq:pi-running-cost-matrix}
    W_{i}(K):= N_{i}+\sum_{j = 1}^{N}K_{j}^{\top} R_{ij}K_{j},\qquad i \in \mathbb N.
\end{equation*}
Throughout, for a smooth map \(f\) between finite-dimensional normed vector spaces, \(Df(x)\) denotes the Fr\'echet derivative of \(f\) at the point \(x\), and \(Df(x)[\Delta x]\) denotes its evaluation at a perturbation \(\Delta x\). For a function of multiple variables, \(\partial_{y} f(x,y)\) stands for the partial Fr\'echet derivative of \(f\) with respect to \(y\), with the remaining parameters held fixed. The subscript \([\,\cdot\,]_{j}\) selects the \(j\)-th block component, and \([\,\cdot\,]_{ij}\) denotes the \((i,j)\)-th block entry of a block-partitioned object.

\begin{definition}
    \label{def:stabilizing-profiles}
    The set of mean-square stabilizing feedback profiles is defined by
    \begin{equation*}
        \label{eq:pi-stabilizing-set}
        \mathcal{G}
        :=
        \left\{
        K = (K_{1},\ldots,K_{N}) \in \mathcal{K}:
        (K_{1},\ldots,K_{N})
        \text{ is a mean-square stabilizer of }
        \eqref{eq:game-dynamics}
        \right\}.
    \end{equation*}
\end{definition}
\begin{lemma}
    \label{lem:pi-policy-evaluation}
    For every \(K \in \mathcal{G}\), the linear operator \(\mathcal{L}_{K}:\mathcal{S}^{n} \to \mathcal{S}^{n}\) is invertible. Moreover, for every \(i \in \mathbb N\), the Lyapunov equation
    \begin{equation*}
        \label{eq:pi-evaluation-equation}
        \mathcal{L}_{K}(P_{i})+W_{i}(K) = 0
    \end{equation*}
    admits a unique solution \(P_{i} \in \mathcal{S}_{++}^{n}\).
\end{lemma}
\begin{proof}
    Fix \(K \in \mathcal{G}\). By Lemma~\ref{lem:mean-square-stability}, for every \(\Lambda \in \mathcal{S}^{n}\), the equation \(\mathcal{L}_{K}(P) = -\Lambda\) admits a unique solution \(P \in \mathcal{S}^{n}\). Therefore, \(\mathcal{L}_{K}\) is bijective and hence invertible. For every \(i \in \mathbb N\), since \(N_{i} \in \mathcal{S}_{++}^{n}\) and \(R_{ij} \in \mathcal{S}_{++}^{m_{j}}\), we have \(W_{i}(K) \in \mathcal{S}_{++}^{n}\). Applying Lemma~\ref{lem:mean-square-stability} with \(\Lambda = W_{i}(K)\) shows that
    \begin{equation*}
        \label{eq:pi-evaluation-solvability}
        \mathcal{L}_{K}(P_{i})+W_{i}(K) = 0
    \end{equation*}
    has a unique solution \(P_{i} \in \mathcal{S}_{++}^{n}\).
\end{proof}
Based on the invertibility of \(\mathcal{L}_{K}\) in Lemma~\ref{lem:pi-policy-evaluation}, for \(K \in \mathcal{G}\), we define the policy evaluation map
\begin{equation*}
    \label{eq:pi-evaluation-map}
    \mathcal{P}_{i}(K):= -\mathcal{L}_{K}^{-1}\bigl(W_{i}(K)\bigr),
    \qquad i \in \mathbb N.
\end{equation*}
For \(K \in \mathcal{G}\), the stationarity condition for player \(i\) with value matrix \(P_{i} = \mathcal{P}_{i}(K)\) reads
\begin{equation}
    \label{eq:pi-stationarity}
    \begin{aligned}
        & \left(R_{ii}+D_{i}^{\top} P_{i}D_{i}+\int_{\mathcal{E}}F_{i}(e)^{\top} P_{i}F_{i}(e)\lambda(\mathrm{d}e)\right)K_{i}\\
        & +\sum_{j \neq i}\left(D_{i}^{\top} P_{i}D_{j}+\int_{\mathcal{E}}F_{i}(e)^{\top} P_{i}F_{j}(e)\lambda(\mathrm{d}e)\right)K_{j}\\
        & =-\left(B_{i}^{\top} P_{i}+D_{i}^{\top} P_{i}C+\int_{\mathcal{E}}F_{i}(e)^{\top} P_{i}E(e)\lambda(\mathrm{d}e)\right).
    \end{aligned}
\end{equation}
Since \(R_{ii} \in \mathcal{S}_{++}^{m_{i}}\) and \(P_{i} \in \mathcal{S}_{++}^{n}\), the diagonal block
\begin{equation*}
    \label{eq:pi-diagonal-positivity}
    R_{ii}+D_{i}^{\top} P_{i}D_{i}+\int_{\mathcal{E}}F_{i}(e)^{\top} P_{i}F_{i}(e)\lambda(\mathrm{d}e)
\end{equation*}
is positive definite and hence invertible.

Stacking~\eqref{eq:pi-stationarity} for \(i = 1,\ldots,N\) gives the linear system
\begin{equation*}
    \label{eq:pi-stacked-stationarity}
    \mathcal{M}(P_{1},\ldots,P_{N})
    \begin{bmatrix}K_{1}\\K_{2}\\\vdots\\K_{N}\end{bmatrix}
    = -\mathcal{N}(P_{1},\ldots,P_{N}),
\end{equation*}
where \(\mathcal{M}(P_{1},\ldots,P_{N})\) is the \(N\times N\) block matrix with blocks
\begin{equation}
    \label{eq:pi-improvement-matrix}
    [\mathcal{M}(P_{1},\ldots,P_{N})]_{ij}:=
    \begin{cases}
        R_{ii}+D_{i}^{\top} P_{i}D_{i}+\int_{\mathcal{E}}F_{i}(e)^{\top} P_{i}F_{i}(e)\lambda(\mathrm{d}e), & j = i,    \\[2mm]
        D_{i}^{\top} P_{i}D_{j}+\int_{\mathcal{E}}F_{i}(e)^{\top} P_{i}F_{j}(e)\lambda(\mathrm{d}e),        & j \neq i,
    \end{cases}
\end{equation}
and \(\mathcal{N}(P_{1},\ldots,P_{N})\) is the stacked matrix
\begin{equation}
    \label{eq:pi-improvement-right-hand-side}
    [\mathcal{N}(P_{1},\ldots,P_{N})]_{i}:=
    B_{i}^{\top} P_{i}+D_{i}^{\top} P_{i}C+\int_{\mathcal{E}}F_{i}(e)^{\top} P_{i}E(e)\lambda(\mathrm{d}e).
\end{equation}
When \(\mathcal{M}(P_{1},\ldots,P_{N})\) is invertible, the policy improvement is
\begin{equation*}
    \label{eq:pi-policy-improvement}
    \begin{bmatrix}K_{1}\\K_{2}\\\vdots\\K_{N}\end{bmatrix}
    = -\mathcal{M}(P_{1},\ldots,P_{N})^{-1}\mathcal{N}(P_{1},\ldots,P_{N}).
\end{equation*}
Define the policy iteration map
\begin{equation}
    \label{eq:pi-iteration-map}
    \mathcal{T}(K):=
    -\mathcal{M}\bigl(\mathcal{P}_{1}(K),\ldots,\mathcal{P}_{N}(K)\bigr)^{-1}
    \mathcal{N}\bigl(\mathcal{P}_{1}(K),\ldots,\mathcal{P}_{N}(K)\bigr),
\end{equation}
where the stacked matrix is identified with its gain profile, and define
\begin{equation}
    \label{eq:pi-admissible-domain}
    \mathcal{D}:=
    \left\{K \in \mathcal{G}:
    \mathcal{M}\bigl(\mathcal{P}_{1}(K),\ldots,\mathcal{P}_{N}(K)\bigr)
    \ \text{is nonsingular}\right\}.
\end{equation}

We now present a partially model-free RL algorithm based on integral policy evaluation and simultaneous policy improvement.

\begin{algorithm}[H]
    \caption{Partially model-free policy-iteration RL for the \(N\)-player Problem (LQ NZSSDG)}
    \label{alg:policy-iteration}
    \textbf{Initialize:} Select \(K^{(0)} = (K_{1}^{(0)},\ldots,K_{N}^{(0)}) \in \mathcal{D}\). Fix \(t \geq 0\) and \(\Delta t > 0\).
    \\[2mm]
    \textbf{For \(l = 0,1,2,\ldots\):}
    \begin{itemize}
        \item \textbf{Data Collection:} Generate the closed-loop trajectory \(X^{(l)}(\cdot)\) under \(u_{j}^{(l)}(s) = K_{j}^{(l)}X^{(l)}(s)\), \(j \in \mathbb N\), and select sampling windows \([t,t+\Delta t]\).
        \item \textbf{Policy Evaluation:} For every \(i \in \mathbb N\), solve for \(P_{i}^{(l+1)}\) using
        \begin{equation}
            \label{eq:pi-integral-evaluation}
            \begin{aligned}
                & X^{(l)}(t)^{\top} P_{i}^{(l+1)}X^{(l)}(t)
                    -\mathbb{E}^{\mathcal{F}_{t}}\!\left[
                        X^{(l)}(t+\Delta t)^{\top} P_{i}^{(l+1)}
                        X^{(l)}(t+\Delta t)\right]\\
                & \qquad = \mathbb{E}^{\mathcal{F}_{t}}\!\int_{t}^{t+\Delta t}
                    X^{(l)}(s)^{\top} W_{i}(K^{(l)})X^{(l)}(s)\,\mathrm{d}s,
                    \qquad i \in \mathbb N.
            \end{aligned}
        \end{equation}
        \item \textbf{Policy Improvement:} Update \((K_{1}^{(l+1)},\ldots,K_{N}^{(l+1)})\) by the formula
        \begin{equation}
        \label{eq:pi-gain-update}
            \begin{bmatrix}K_{1}^{(l+1)}\\\vdots\\K_{N}^{(l+1)}\end{bmatrix}
            = -\mathcal{M}(P_{1}^{(l+1)},\ldots,P_{N}^{(l+1)})^{-1}\mathcal{N}(P_{1}^{(l+1)},\ldots,P_{N}^{(l+1)}).
        \end{equation}
    \end{itemize}
    \textbf{Stop when} \(\sum_{i = 1}^{N} \|P_{i}^{(l+1)}-P_{i}^{(l)}\|_F < \varepsilon\) for \(l \geq 1\).
\end{algorithm}

For the convenience of subsequent analysis, denote
\begin{equation*}
    \label{eq:pi-iterate-coefficients}
    A_{l}:= A+\sum_{j = 1}^{N}B_{j}K_{j}^{(l)},\quad
    C_{l}:= C+\sum_{j = 1}^{N}D_{j}K_{j}^{(l)},\quad
    E_{l}(e):= E(e)+\sum_{j = 1}^{N}F_{j}(e)K_{j}^{(l)},\quad l \geq 0.
\end{equation*}

\begin{lemma}[Equivalence of policy evaluation and the Lyapunov recursion]
    \label{lem:pi-evaluation-equivalence}
    For every \(i \in \mathbb N\), solving~\eqref{eq:pi-integral-evaluation} for all initial states and all sufficiently small \(\Delta t > 0\) is equivalent to solving
    \begin{equation}
        \label{eq:pi-lyapunov-recursion}
        \mathcal{L}_{K^{(l)}}(P_{i}^{(l+1)})+W_{i}(K^{(l)}) = 0.
    \end{equation}
\end{lemma}

\begin{proof}
    Under the feedback control \(u_{j}^{(l)}(s) = K_{j}^{(l)}X^{(l)}(s)\), \(j \in \mathbb N\), the state satisfies
    \begin{equation*}
        \label{eq:pi-iterate-dynamics}
        \left\{
        \begin{aligned}
            \mathrm{d}X^{(l)}(s) & =A_{l}X^{(l)}(s)\,\mathrm{d}s+C_{l}X^{(l)}(s)\,\mathrm{d}W(s)
            +\int_{\mathcal{E}}E_{l}(e)X^{(l)}(s)\widetilde{N}(\mathrm{d}s,\mathrm{d}e),\quad s \geq t,\\
            X^{(l)}(t) & =x.
        \end{aligned}
        \right.
    \end{equation*}
    Applying the It\^o formula to \(X^{(l)}(s)^{\top} P_{i}^{(l+1)}X^{(l)}(s)\) yields
    \begin{equation}
        \label{eq:pi-iterate-ito}
        \begin{aligned}
            & \mathrm{d}\!\Bigl[X^{(l)}(s)^{\top} P_{i}^{(l+1)}X^{(l)}(s)\Bigr] \\
            & \quad = X^{(l)}(s)^{\top} \Bigl\{A_{l}^{\top} P_{i}^{(l+1)}+P_{i}^{(l+1)}A_{l} +C_{l}^{\top} P_{i}^{(l+1)}C_{l}
            +\int_{\mathcal{E}}E_{l}(e)^{\top} P_{i}^{(l+1)}E_{l}(e)\lambda(\mathrm{d}e)\Bigr\} X^{(l)}(s)\,\mathrm{d}s\\
            & \qquad+\{\cdots\}\,\mathrm{d}W(s) + \{\cdots\}\,\widetilde{N}(\mathrm{d}s,\mathrm{d}e).
        \end{aligned}
    \end{equation}
    Integrating~\eqref{eq:pi-iterate-ito} from \(t\) to \(t+\Delta t\) gives
    \begin{equation}
        \label{eq:pi-integrated-ito}
        \begin{aligned}
            & X^{(l)} (t+\Delta t)^{\top} P_{i}^{(l+1)} X^{(l)}(t+\Delta t)-x^{\top} P_{i}^{(l+1)}x \\
            & \quad = \int_{t}^{t+\Delta t} X^{(l)}(s)^{\top} \Bigl(A_{l}^{\top} P_{i}^{(l+1)}+P_{i}^{(l+1)}A_{l}+C_{l}^{\top} P_{i}^{(l+1)}C_{l}
            +\int_{\mathcal{E}}E_{l}(e)^{\top} P_{i}^{(l+1)}E_{l}(e) \lambda(\mathrm{d}e)\Bigr) X^{(l)}(s)\,\mathrm{d}s\\
            & \qquad+\int_{t}^{t+\Delta t}\{\cdots\}\,\mathrm{d}W(s)+\int_{t}^{t+\Delta t}\int_{\mathcal{E}}\{\cdots\}\widetilde{N}(\mathrm{d}s,\mathrm{d}e).
        \end{aligned}
    \end{equation}
    Taking the conditional expectation \(\mathbb{E}^{\mathcal{F}_{t}}\) on both sides of~\eqref{eq:pi-integrated-ito}, we obtain
    \begin{equation}
        \label{eq:pi-conditional-expectation}
        \begin{aligned}
            & x^{\top} P_{i}^{(l+1)}x-\mathbb{E}^{\mathcal{F}_{t}}\!\Bigl[X^{(l)}(t+\Delta t)^{\top} P_{i}^{(l+1)}X^{(l)}(t+\Delta t)\Bigr] \\
            & \quad = -\mathbb{E}^{\mathcal{F}_{t}}\!\int_{t}^{t+\Delta t} X^{(l)}(s)^{\top} \Bigl(A_{l}^{\top} P_{i}^{(l+1)}+P_{i}^{(l+1)}A_{l}+C_{l}^{\top} P_{i}^{(l+1)}C_{l}
            +\int_{\mathcal{E}}E_{l}(e)^{\top} P_{i}^{(l+1)}E_{l}(e) \lambda(\mathrm{d}e)\Bigr) X^{(l)}(s)\,\mathrm{d}s.
        \end{aligned}
    \end{equation}
    The Lyapunov recursion~\eqref{eq:pi-lyapunov-recursion} is exactly \(A_{l}^{\top} P_{i}^{(l+1)}+P_{i}^{(l+1)}A_{l}+C_{l}^{\top} P_{i}^{(l+1)}C_{l}+\int_{\mathcal{E}} E_{l}(e)^{\top} P_{i}^{(l+1)}E_{l}(e) \lambda(\mathrm{d}e) = -W_{i}(K^{(l)})\). Substitution of this equality into~\eqref{eq:pi-conditional-expectation} gives~\eqref{eq:pi-integral-evaluation}.

    Conversely, suppose that \(P_{i}^{(l+1)}\) satisfies the policy evaluation equation~\eqref{eq:pi-integral-evaluation} for all \(x\) and for any \(\Delta t > 0\). Comparing~\eqref{eq:pi-integral-evaluation} with~\eqref{eq:pi-conditional-expectation}, we obtain
    \begin{equation*}
        \label{eq:pi-converse-integral}
        \begin{aligned}
            \mathbb{E}^{\mathcal{F}_{t}}
            \left[\int_{t}^{t+\Delta t} X^{(l)}(s)^{\top}
                W_{i}(K^{(l)}) X^{(l)}(s)\,\mathrm{d}s\right]
            & +\mathbb{E}^{\mathcal{F}_{t}}\!\int_{t}^{t+\Delta t} X^{(l)}(s)^{\top} \left(A_{l}^{\top} P_{i}^{(l+1)}+P_{i}^{(l+1)}A_{l} \right.\\
            & \quad \left. +C_{l}^{\top} P_{i}^{(l+1)}C_{l} +\int_{\mathcal{E}}E_{l}(e)^{\top} P_{i}^{(l+1)}E_{l}(e) \lambda(\mathrm{d}e)\right)
            X^{(l)}(s)\,\mathrm{d}s = 0.
        \end{aligned}
    \end{equation*}
    Dividing both sides of the above identity by \(\Delta t\), we obtain
    \begin{equation*}
        \label{eq:pi-converse-average}
        \begin{aligned}
            \frac{1}{\Delta t}\mathbb{E}^{\mathcal{F}_{t}}\!\int_{t}^{t+\Delta t} X^{(l)}(s)^{\top}
            \Big\{W_{i}(K^{(l)}) & +A_{l}^{\top} P_{i}^{(l+1)}+P_{i}^{(l+1)}A_{l}+C_{l}^{\top} P_{i}^{(l+1)}C_{l} \\
            & \quad +\int_{\mathcal{E}}E_{l}(e)^{\top} P_{i}^{(l+1)}E_{l}(e) \lambda(\mathrm{d}e)\Big\}
            X^{(l)}(s)\,\mathrm{d}s = 0.
        \end{aligned}
    \end{equation*}
    Letting \(\Delta t \rightarrow 0\), we obtain
    \begin{equation*}
        \label{eq:pi-vanishing-quadratic-form}
        x^{\top}\Big(A_{l}^{\top} P_{i}^{(l+1)}+P_{i}^{(l+1)}A_{l}+C_{l}^{\top} P_{i}^{(l+1)}C_{l}
        +\int_{\mathcal{E}}E_{l}(e)^{\top} P_{i}^{(l+1)}E_{l}(e) \lambda(\mathrm{d}e)+W_{i}(K^{(l)})\Big)x = 0,
        \quad \forall x \in \mathbb{R}^{n}.
    \end{equation*}
    Therefore, \(A_{l}^{\top} P_{i}^{(l+1)}+P_{i}^{(l+1)}A_{l}+C_{l}^{\top} P_{i}^{(l+1)}C_{l} +\int_{\mathcal{E}}E_{l}(e)^{\top} P_{i}^{(l+1)}E_{l}(e) \lambda(\mathrm{d}e)+W_{i}(K^{(l)}) = 0\), which is exactly~\eqref{eq:pi-lyapunov-recursion}. The proof is complete.
\end{proof}

By Lemma~\ref{lem:pi-evaluation-equivalence}, whenever \(K^{(l)} \in \mathcal{D}\), the policy evaluation and policy improvement steps in Algorithm~\ref{alg:policy-iteration} can be written as \(P_{i}^{(l+1)} = \mathcal{P}_{i}(K^{(l)}),\ i \in \mathbb N,\ K^{(l+1)} = \mathcal{T}(K^{(l)})\).

\begin{lemma}[Openness and smoothness]
    \label{lem:pi-smoothness}
    Both \(\mathcal{G}\) and \(\mathcal{D}\) are open subsets of \(\mathcal{K}\). For each \(i\in\mathbb N\), \(\mathcal{P}_{i}\colon \mathcal{G}\to \mathcal{S}_{++}^{n}\) is \(C^{\infty}\), and the mapping \(\mathcal{T}\colon\mathcal{D}\to\mathcal{K}\) is also \(C^{\infty}\), i.e., infinitely continuously differentiable.
\end{lemma}
\begin{proof}
    Take any \(K \in \mathcal{G}\). By~\cite[Theorem 1.5.3]{damm2004rational}, mean-square stability is equivalent to the existence of some \(X > 0\) satisfying the strict Lyapunov inequality \(\mathcal{L}_{K}(X) < 0\). Since the map \(K\mapsto \mathcal{L}_{K}\) is polynomial, hence continuous, the strict inequality \(\mathcal{L}_{K}(X) < 0\) remains valid under sufficiently small perturbations of \(K\). Thus \(\mathcal{G}\) is open in \(\mathcal{K}\).

    Define \(\Phi_{i}(K,P):= \mathcal{L}_{K}(P)+W_{i}(K)\) for \((K,P) \in \mathcal{K}\times\mathcal{S}^{n}\). The mapping \(\Phi_{i}\) is polynomial in the pair \((K,P)\), hence \(C^{\infty}\). Because \(\mathcal{L}_{K}\) is linear in \(P\) and \(W_{i}(K)\) is independent of \(P\), the partial Fréchet derivative with respect to \(P\) satisfies \(\partial_{P}\Phi_{i}(K,P) = \mathcal{L}_{K}\). By Lemma~\ref{lem:pi-policy-evaluation}, \(\mathcal{L}_{K}\) is invertible for every \(K \in \mathcal{G}\). For each \(K \in \mathcal{G}\), the equation \(\Phi_{i}(K,\mathcal{P}_{i}(K)) = 0\) has a unique solution. Applying the implicit function theorem, we conclude that \(\mathcal{P}_{i}:\mathcal{G} \to \mathcal{S}^{n}\) is \(C^{\infty}\). By uniqueness of the policy-evaluation solution, this locally constructed solution coincides globally with \(\mathcal{P}_{i}(K)\). Therefore, \(\mathcal{P}_{i} \in C^{\infty}(\mathcal{G};\mathcal{S}_{++}^{n})\).

    The openness of \(\mathcal{G}\) and continuity of every \(\mathcal{P}_{i}\) have just been established. Since the set of nonsingular matrices is open, the defining condition in~\eqref{eq:pi-admissible-domain} shows that \(\mathcal{D}\) is open. The maps \(\mathcal{M}\) and \(\mathcal{N}\) are polynomial in their matrix arguments, while matrix inversion is smooth on the set of nonsingular matrices. Formula~\eqref{eq:pi-iteration-map} therefore proves that \(\mathcal{T}\) is \(C^{\infty}\) on \(\mathcal{D}\).
\end{proof}
For \(i,j \in \mathbb N\), denote
\begin{equation}
    \label{eq:pi-stationarity-residual}
    \Xi_{ij}(K):=
    B_{j}^{\top}\mathcal{P}_{i}(K)+D_{j}^{\top}\mathcal{P}_{i}(K)C_{K}
    +\int_{\mathcal{E}}F_{j}(e)^{\top}\mathcal{P}_{i}(K)E_{K}(e)\lambda(\mathrm{d}e)
    +R_{ij}K_{j}.
\end{equation}

\begin{lemma}
    \label{lem:pi-derivatives}
    For \(K \in \mathcal{G}\) and \(\Delta = (\Delta_{1},\ldots,\Delta_{N}) \in \mathcal{K}\),
    \begin{equation}
        \label{eq:pi-evaluation-derivative}
        D\mathcal{P}_{i}(K)[\Delta]
        = -\mathcal{L}_{K}^{-1}\left[
            \sum_{j = 1}^{N}\operatorname{sym}
            \bigl(\Delta_{j}^{\top}\Xi_{ij}(K)\bigr)\right],
        \qquad i \in \mathbb N.
    \end{equation}
    Suppose that \(K^{*} \in \mathcal{D}\) is a fixed point of \(\mathcal{T}\), i.e., \(K^{*} = \mathcal{T}(K^{*})\), and let \(P_{i}^{*} = \mathcal{P}_{i}(K^{*}),\ i \in \mathbb N\). Then, for every \(\Delta = (\Delta_{1},\ldots,\Delta_{N}) \in \mathcal{K}\),
    \begin{equation}
        \label{eq:pi-fixed-point-evaluation-derivative}
        D\mathcal{P}_{i}(K^{*})[\Delta]
        =
        -\mathcal{L}_{K^{*}}^{-1}
        \left[
            \sum_{\substack{j = 1 \\ j \neq i}}^{N}
            \operatorname{sym}
            \bigl(\Delta_{j}^{\top}\Xi_{ij}(K^{*})\bigr)
            \right],
        \qquad i \in \mathbb N.
    \end{equation}
    Moreover,
    \begin{equation}
        \label{eq:pi-fixed-point-map-derivative}
        D\mathcal{T}(K^{*})[\Delta]
        =
        -\mathcal{M}(P_{1}^{*},\ldots,P_{N}^{*})^{-1}
        \begin{bmatrix}
            \Theta_{1}^{*}(\Delta) \\
            \vdots                 \\
            \Theta_{N}^{*}(\Delta)
        \end{bmatrix},
    \end{equation}
    where
    \begin{equation*}
        \label{eq:pi-derivative-right-hand-side}
        \Theta_{i}^{*}(\Delta)
        =
        B_{i}^{\top} D\mathcal{P}_{i}(K^{*})[\Delta]
        +
        D_{i}^{\top}
        D\mathcal{P}_{i}(K^{*})[\Delta]C_{K}^{*}
        +
        \int_{\mathcal{E}}
        F_{i}(e)^{\top}
        D\mathcal{P}_{i}(K^{*})[\Delta]
        E_{K}^{*}(e)\lambda(\mathrm{d}e),
        \qquad i \in \mathbb N.
    \end{equation*}
\end{lemma}

\begin{proof}
    Since \(\mathcal{G}\) is open, for all sufficiently small \(\alpha\) one may set \(K(\alpha) = K+\alpha\Delta \in \mathcal{G}\). Put
    \begin{equation*}
        \label{eq:pi-coefficient-perturbations}
        \Delta_{A} = \sum_{j = 1}^{N}B_{j}\Delta_{j},\qquad
        \Delta_{C} = \sum_{j = 1}^{N}D_{j}\Delta_{j},\qquad
        \Delta_{E}(e) = \sum_{j = 1}^{N}F_{j}(e)\Delta_{j}.
    \end{equation*}
    Then \(A_{K(\alpha)} = A_{K}+\alpha\Delta_{A}\), \(C_{K(\alpha)} = C_{K}+\alpha\Delta_{C}\), and \(E_{K(\alpha)}(e) = E_{K}(e)+\alpha\Delta_{E}(e)\). The jump term therefore satisfies
    \begin{equation*}
        \label{eq:pi-jump-derivative}
        \begin{aligned}
            \left.\frac{d}{d\alpha}\int_{\mathcal{E}}E_{K(\alpha)}(e)^{\top} PE_{K(\alpha)}(e)\lambda(\mathrm{d}e)\right|_{\alpha = 0} = \int_{\mathcal{E}}\left[\Delta_{E}(e)^{\top} PE_{K}(e)+E_{K}(e)^{\top} P\Delta_{E}(e)\right]\lambda(\mathrm{d}e).
        \end{aligned}
    \end{equation*}
    Hence, for \(P \in \mathcal{S}^{n}\),
    \begin{equation}
        \label{eq:pi-lyapunov-derivative}
        \left.\frac{\mathrm{d}}{\mathrm{d}\alpha}
        \mathcal{L}_{K(\alpha)}(P)
        \right|_{\alpha = 0}
        =
        \Delta_{A}^{\top} P+P\Delta_{A}
        +\Delta_{C}^{\top} PC_{K}+C_{K}^{\top} P\Delta_{C}
        +\int_{\mathcal{E}}(\Delta_{E}(e)^{\top} PE_{K}(e)+E_{K}(e)^{\top} P\Delta_{E}(e))\lambda(\mathrm{d}e).
    \end{equation}
    \begin{equation}
        \label{eq:pi-running-cost-derivative}
        \left.\frac{\mathrm{d}}{\mathrm{d}\alpha}
        W_{i}(K(\alpha))
        \right|_{\alpha = 0}
        =
        \sum_{j = 1}^{N}
        \bigl(\Delta_{j}^{\top} R_{ij}K_{j}+K_{j}^{\top} R_{ij}\Delta_{j}\bigr).
    \end{equation}
    Differentiating \(\mathcal{L}_{K(\alpha)}[\mathcal{P}_{i}(K(\alpha))]+W_{i}(K(\alpha)) = 0\) at \(\alpha = 0\) and substituting~\eqref{eq:pi-lyapunov-derivative}--\eqref{eq:pi-running-cost-derivative} give
    \begin{equation*}
        \label{eq:pi-differentiated-evaluation}
        \begin{aligned}
            0 = {} & \mathcal{L}_{K}\bigl(D\mathcal{P}_{i}(K)[\Delta]\bigr)
            +\Delta_{A}^{\top} \mathcal{P}_{i}(K)+\mathcal{P}_{i}(K)\Delta_{A}
            +\Delta_{C}^{\top} \mathcal{P}_{i}(K)C_{K}+C_{K}^{\top} \mathcal{P}_{i}(K)\Delta_{C}\\
            & +\int_{\mathcal{E}}\!\left[
                \Delta_{E}(e)^{\top} \mathcal{P}_{i}(K)E_{K}(e)
                +E_{K}(e)^{\top} \mathcal{P}_{i}(K)\Delta_{E}(e)\right]\lambda(\mathrm{d}e)\\
            & +\sum_{j = 1}^{N}\left(
            \Delta_{j}^{\top} R_{ij}K_{j}+K_{j}^{\top} R_{ij}\Delta_{j}\right).
        \end{aligned}
    \end{equation*}
    Substituting the expressions for \(\Delta_{A}\), \(\Delta_{C}\), and \(\Delta_{E}\) in terms of \(\Delta_{j}\), collecting like terms, and using~\eqref{eq:pi-stationarity-residual}, we obtain
    \begin{equation}
        \label{eq:pi-derivative-lyapunov-equation}
        \mathcal{L}_{K}(D\mathcal{P}_{i}(K)[\Delta]) + \sum_{j = 1}^{N}
        \operatorname{sym}\bigl(\Delta_{j}^{\top}\Xi_{ij}(K)\bigr) = 0.
    \end{equation}
    Since \(\mathcal{L}_{K}\) is invertible, we can derive~\eqref{eq:pi-evaluation-derivative} from~\eqref{eq:pi-derivative-lyapunov-equation}.

    At a fixed point, the derivative formulas simplify as follows. Since \(K^{*} = \mathcal{T}(K^{*})\), the stationarity condition gives \(\Xi_{ii}(K^{*}) = 0,\ i \in \mathbb N\). Thus,~\eqref{eq:pi-fixed-point-evaluation-derivative} follows directly from~\eqref{eq:pi-evaluation-derivative}. Differentiating the defining identity \(\mathcal{M}(\mathcal{P}_{1}(K),\ldots,\mathcal{P}_{N}(K))\mathcal{T}(K)+\mathcal{N}(\mathcal{P}_{1}(K),\ldots,\mathcal{P}_{N}(K))=0\) at \(K^{*}\) along \(\Delta\) gives, for each \(i \in \mathbb N\),
    \begin{equation}
        \label{eq:pi-component-map-derivative}
        \begin{aligned}
            & \left(
            R_{ii}+D_{i}^{\top} P_{i}^{*}D_{i}
            +\int_{\mathcal{E}}
            F_{i}(e)^{\top} P_{i}^{*}F_{i}(e)\lambda(\mathrm{d}e)
            \right)
            [D\mathcal{T}(K^{*})[\Delta]]_{i}\\
            & \quad+
            \sum_{j \neq i}
            \left(
            D_{i}^{\top} P_{i}^{*}D_{j}
            +\int_{\mathcal{E}}
            F_{i}(e)^{\top} P_{i}^{*}F_{j}(e)\lambda(\mathrm{d}e)
            \right)
            [D\mathcal{T}(K^{*})[\Delta]]_{j}\\
            & =
            -B_{i}^{\top} D\mathcal{P}_{i}(K^{*})[\Delta]
            -D_{i}^{\top} D\mathcal{P}_{i}(K^{*})[\Delta]C_{K}^{*}
            -\int_{\mathcal{E}}
            F_{i}(e)^{\top} D\mathcal{P}_{i}(K^{*})[\Delta]
            E_{K}^{*}(e)\lambda(\mathrm{d}e).
        \end{aligned}
    \end{equation}
    Stacking~\eqref{eq:pi-component-map-derivative} for \(i = 1,\ldots,N\) yields
    \begin{equation*}
        \label{eq:pi-stacked-map-derivative}
        \mathcal{M}(P_{1}^{*},\ldots,P_{N}^{*})
        D\mathcal{T}(K^{*})[\Delta]
        =
        -
        \begin{bmatrix}
            \Theta_{1}^{*}(\Delta) \\
            \vdots                 \\
            \Theta_{N}^{*}(\Delta)
        \end{bmatrix}.
    \end{equation*}
    Since \(K^{*} \in \mathcal{D}\), \(\mathcal{M}(P_{1}^{*},\ldots,P_{N}^{*})\) is nonsingular. Therefore,~\eqref{eq:pi-fixed-point-map-derivative} follows.
\end{proof}

The traditional fixed-point-centered sufficient condition is therefore that the spectral radius of \(D\mathcal{T}(K^{*})\) is less than one, i.e., \(\rho(D\mathcal{T}(K^{*})) < 1\). The theorem below instead uses a derivative bound on a ball centered at the selected \(K^{(0)}\), so neither \(K^{*}\) nor this spectral condition needs to be assumed in advance.

\begin{assumption}
    \label{ass:invariant-ball}
    There exist an initial feedback profile \(K^{(0)} \in \mathcal{D}\), a constant \(\delta > 0\), a constant \(q \in [0,1)\), and a norm \(\|\cdot\|_{*}\) on \(\mathcal{K}\) such that
    \begin{enumerate}
        \item[(a)] \(\bar{B}:= \{K \in \mathcal{K}:\|K-K^{(0)}\|_{*} \leq \delta\} \subset \mathcal{D}\);
        \item[(b)] \(\sup_{K \in \bar{B}} \|D\mathcal{T}(K)\|_{*} \leq q\);
        \item[(c)] \(\eta:= \|\mathcal{T}(K^{(0)})-K^{(0)}\|_{*} \leq (1-q)\delta\).
    \end{enumerate}
\end{assumption}

We now present the following theorem, which establishes local geometric convergence of the feedback gains and value matrices and identifies the limiting profile as a closed-loop Nash equilibrium.

\begin{theorem}
    \label{thm:pi-convergence}
    Suppose Assumption~\ref{ass:invariant-ball} holds. Then:
    \begin{enumerate}
        \item[(i)] \(\mathcal{T}\) has a unique fixed point \(K^{*} \in \bar{B}\), and \(\|K^{*}-K^{(0)}\|_{*} \leq \eta/(1-q)\).
        \item[(ii)] Every iterate of Algorithm~\ref{alg:policy-iteration} is well defined and mean-square stabilizing, with
            \begin{equation*}
                \label{eq:pi-gain-convergence-bound}
                \|K^{(l)}-K^{*}\|_{*}
                \leq \frac{q^{l}\eta}{1-q},\qquad l \geq 0.
            \end{equation*}
        \item[(iii)] For each \(i \in \mathbb N\), with \(P_{i}^{*} = \mathcal{P}_{i}(K^{*})\), there exists \(L_{i} > 0\) such that
            \begin{equation*}
                \label{eq:pi-value-convergence-bound}
                \|P_{i}^{(l+1)}-P_{i}^{*}\|_F
                \leq L_{i}\frac{q^{l}\eta}{1-q},
                \qquad P_{i}^{(l+1)} \in \mathcal{S}_{++}^{n}.
            \end{equation*}
        \item[(iv)] The feedback profile \((K_{1}^{*}X,\ldots,K_{N}^{*}X)\) is a closed-loop Nash equilibrium of Problem (LQ NZSSDG).
    \end{enumerate}
\end{theorem}

\begin{proof}
    \begin{enumerate}
        \item[(i)] By Assumption~\ref{ass:invariant-ball}(b),
            \begin{equation*}
                \label{eq:pi-uniform-derivative-bound}
                \|D\mathcal{T}(K)\|_{*}
                \leq q < 1,
                \qquad K \in \bar{B}.
            \end{equation*}
            Since \(\bar{B}\) is convex, the mean-value inequality gives, for any \(K,K' \in \bar{B}\),
            \begin{equation}
                \label{eq:pi-contraction}
                \|\mathcal{T}(K)-\mathcal{T}(K')\|_{*}
                \leq
                q\|K-K'\|_{*}.
            \end{equation}
            Hence \(\mathcal{T}\) is a \(q\)-contraction on \(\bar{B}\). For any \(K \in \bar{B}\), the triangle inequality gives
            \begin{equation*}
                \label{eq:pi-triangle-inequality}
                \|\mathcal{T}(K)-K^{(0)}\|_{*}
                \leq \|\mathcal{T}(K)-\mathcal{T}(K^{(0)})\|_{*}
                +\|\mathcal{T}(K^{(0)})-K^{(0)}\|_{*}.
            \end{equation*}
            From~\eqref{eq:pi-contraction} and Assumption~\ref{ass:invariant-ball}(c),
            \begin{equation*}
                \label{eq:pi-map-displacement-bound}
                \|\mathcal{T}(K)-\mathcal{T}(K^{(0)})\|_{*}
                \leq q\|K-K^{(0)}\|_{*} \leq q\delta
            \end{equation*}
            and
            \begin{equation*}
                \label{eq:pi-initial-residual-bound}
                \|\mathcal{T}(K^{(0)})-K^{(0)}\|_{*}
                = \eta \leq (1-q)\delta.
            \end{equation*}
            Consequently,
            \begin{equation*}
                \label{eq:pi-ball-displacement-bound}
                \|\mathcal{T}(K)-K^{(0)}\|_{*}
                \leq q\delta+(1-q)\delta = \delta,
            \end{equation*}
            so \(\mathcal{T}(K) \in \bar{B}\) and
            \begin{equation}
                \label{eq:pi-ball-invariance}
                \mathcal{T}(\bar{B}) \subset \bar{B}.
            \end{equation}

            Since \(\mathcal{K}\) is finite dimensional, the closed ball \(\bar{B}\) is complete. Hence, by~\eqref{eq:pi-contraction},~\eqref{eq:pi-ball-invariance}, and Banach's fixed-point theorem, \(\mathcal{T}\) admits a unique fixed point \(K^{*} \in \bar{B}\). Moreover, since \(K^{*} = \mathcal{T}(K^{*})\),
            \begin{equation*}
                \label{eq:pi-fixed-point-distance-estimate}
                \begin{aligned}
                    \|K^{*}-K^{(0)}\|_{*}
                    & \leq
                    \|\mathcal{T}(K^{*})-\mathcal{T}(K^{(0)})\|_{*}
                    +
                    \|\mathcal{T}(K^{(0)})-K^{(0)}\|_{*}
                    \\
                    & \leq
                    q\|K^{*}-K^{(0)}\|_{*}+\eta.
                \end{aligned}
            \end{equation*}
            Thus,
            \begin{equation*}
                \label{eq:pi-fixed-point-distance-bound}
                \|K^{*}-K^{(0)}\|_{*}
                \leq
                \frac{\eta}{1-q}.
            \end{equation*}
            This proves part~(i).

        \item[(ii)] By Assumption~\ref{ass:invariant-ball}, \(K^{(0)} \in \bar{B}\), and~\eqref{eq:pi-ball-invariance} implies recursively that \(K^{(l)} \in \bar{B},\ l \geq 0\). Assumption~\ref{ass:invariant-ball}(a) gives \(\bar{B} \subset \mathcal{D}\), and by definition \(\mathcal{D} \subset \mathcal{G}\). Hence every \(K^{(l)} \in \bar{B}\) is mean-square stabilizing. Membership in \(\mathcal{D}\) also guarantees that \(\mathcal{M}\bigl(\mathcal{P}_{1}(K^{(l)}),\ldots, \mathcal{P}_{N}(K^{(l)})\bigr)\) is nonsingular, so each policy-improvement step is well defined. Furthermore, using \(K^{(l+1)} = \mathcal{T}(K^{(l)})\) and~\eqref{eq:pi-contraction}, we obtain
            \begin{equation*}
                \label{eq:pi-successive-increment-bound}
                \begin{aligned}
                    \|K^{(l+1)}-K^{(l)}\|_{*}
                    & =
                    \|\mathcal{T}(K^{(l)})
                    -\mathcal{T}(K^{(l-1)})\|_{*}
                    \\
                    & \leq
                    q\|K^{(l)}-K^{(l-1)}\|_{*}\\
                    & \leq
                    q^{l}
                    \|K^{(1)}-K^{(0)}\|_{*}
                    =
                    q^{l}\eta.
                \end{aligned}
            \end{equation*}
            For any \(r > l\),
            \begin{equation*}
                \label{eq:pi-cauchy-estimate}
                \|K^{(r)}-K^{(l)}\|_{*}
                \leq \sum_{j = l}^{r-1}\|K^{(j+1)}-K^{(j)}\|_{*}
                \leq \eta\sum_{j = l}^{r-1}q^{j}.
            \end{equation*}
            Letting \(r \to \infty\) and using \(K^{(r)} \to K^{*}\) yields
            \begin{equation}
                \label{eq:pi-geometric-gain-bound}
                \|K^{(l)}-K^{*}\|_{*}
                \leq
                \frac{q^{l}\eta}{1-q},
                \qquad l \geq 0.
            \end{equation}
            This proves part~(ii).

        \item[(iii)] By Lemma~\ref{lem:pi-smoothness}, each mapping \(\mathcal{P}_{i}\) is continuously differentiable on \(\mathcal{D}\). Since \(\bar{B}\) is compact and convex and \(\bar{B} \subset \mathcal{D}\), there exists a constant \(L_{i} > 0\) such that
            \begin{equation}
                \label{eq:pi-evaluation-lipschitz-bound}
                \|\mathcal{P}_{i}(K)-\mathcal{P}_{i}(K')\|_F
                \leq
                L_{i}\|K-K'\|_{*},
                \qquad K,K' \in \bar{B}.
            \end{equation}
            Since \(P_{i}^{(l+1)} = \mathcal{P}_{i}(K^{(l)})\) and \(P_{i}^{*} = \mathcal{P}_{i}(K^{*})\), equations~\eqref{eq:pi-geometric-gain-bound}--\eqref{eq:pi-evaluation-lipschitz-bound} imply
            \begin{equation*}
                \label{eq:pi-value-error-estimate}
                \|P_{i}^{(l+1)}-P_{i}^{*}\|_F
                \leq L_{i}\|K^{(l)}-K^{*}\|_{*}
                \leq L_{i}\frac{q^{l}\eta}{1-q}.
            \end{equation*}
            Moreover, since \(K^{(l)} \in \bar{B} \subset \mathcal{D}\), Lemma~\ref{lem:pi-policy-evaluation} yields \(P_{i}^{(l+1)} \in \mathcal{S}_{++}^{n}\). This proves part~(iii).

        \item[(iv)] At the fixed point, the evaluation relation is
            \begin{equation}
                \label{eq:pi-fixed-point-lyapunov}
                \mathcal{L}_{K^{*}}(P_{i}^{*})+W_{i}(K^{*}) = 0,
                \qquad i \in \mathbb N.
            \end{equation}
            On the other hand, the fixed-point relation \(K^{*} = \mathcal{T}(K^{*})\) implies the stationarity conditions
            \begin{equation}
                \label{eq:pi-fixed-point-stationarity}
                \begin{aligned}
                    0 = {} &
                    B_{i}^{\top} P_{i}^{*}
                    +
                    D_{i}^{\top} P_{i}^{*}
                    \left(
                    C+\sum_{j = 1}^{N}D_{j}K_{j}^{*}
                    \right)
                    \\
                    & +
                    \int_{\mathcal{E}}
                    F_{i}(e)^{\top} P_{i}^{*}
                    \left(
                    E(e)+
                    \sum_{j = 1}^{N}F_{j}(e)K_{j}^{*}
                    \right)
                    \lambda(\mathrm{d}e)
                    +
                    R_{ii}K_{i}^{*},
                    \qquad i \in \mathbb N.
                \end{aligned}
            \end{equation}

            Furthermore,
            \begin{equation*}
                \label{eq:pi-fixed-point-convexity}
                R_{ii}+D_{i}^{\top} P_{i}^{*}D_{i}
                +\int_{\mathcal{E}}F_{i}(e)^{\top} P_{i}^{*}F_{i}(e)\lambda(\mathrm{d}e)
                > 0
            \end{equation*}
            because \(R_{ii} \in \mathcal{S}_{++}^{m_{i}}\) and \(P_{i}^{*} \in \mathcal{S}_{++}^{n}\). Equations~\eqref{eq:pi-fixed-point-lyapunov}--\eqref{eq:pi-fixed-point-stationarity} therefore satisfy the equilibrium characterization in Lemma~\ref{lem:nash-characterization}. Thus \((K_{1}^{*}X,\ldots,K_{N}^{*}X)\) is a closed-loop Nash equilibrium of Problem~(LQ NZSSDG). This proves part~(iv).
    \end{enumerate}
\end{proof}

\begin{remark}
    \label{rem:pi-local-convergence}
    Theorem~\ref{thm:pi-convergence} is a local convergence theorem, which assumes the existence of a ball centered at a selected \(K^{(0)}\). Within this ball, Algorithm~\ref{alg:policy-iteration} is convergent and can ultimately reach a Nash equilibrium. The location of \(K^{*}\) is not needed; its existence, and the fact that it is a Nash equilibrium, are conclusions.
\end{remark}

\begin{remark}
    \label{rem:pi-stability-invariance}
    Assumption~\ref{ass:invariant-ball}(a) places the entire ball inside the mean-square stabilizing set \(\mathcal{G}\). Combined with the invariance \(\mathcal{T}(\bar{B}) \subset \bar{B}\) established in the proof, this guarantees that every iterate is stabilizing. Thus stability follows from domain inclusion and invariance, whereas geometric convergence follows from the separate derivative bound in Assumption~\ref{ass:invariant-ball}(b).
\end{remark}

Assumption~\ref{ass:invariant-ball} requires a stabilizing ball on which the policy-improvement matrix is nonsingular and the iteration map is contractive. Openness of \(\mathcal D\) alone provides neither a usable radius nor a contraction factor. We therefore give sufficient tests for these requirements and develop a verifiable local convergence certificate.

Fix a candidate \(K^{(0)}\in\mathcal D\) and a radius \(\delta>0\). In this certificate, take the norm in Assumption~\ref{ass:invariant-ball} to be the product spectral norm
\begin{equation*}
    \label{eq:cert-product-norm}
    \|K\|_*:=\Bigl(\sum_{j=1}^{N}\|K_j\|_2^2\Bigr)^{1/2}.
\end{equation*}
Write \(C_0:=C_{K^{(0)}}\) and \(E_0(e):=E_{K^{(0)}}(e)\). Assume that there exist \(S\in\mathcal S_{++}^n\) and \(\alpha>0\) such that
\begin{equation}
    \label{eq:cert-center-lyapunov}
    \mathcal L_{K^{(0)}}(S)\leq-\alpha I_n.
\end{equation}
Define
\begin{equation*}
    \label{eq:cert-perturbation-constants}
    \begin{aligned}
        \beta_B & :=\Bigl(\sum_j\|B_j\|_2^2\Bigr)^{1/2},
        & \beta_D & :=\Bigl(\sum_j\|D_j\|_2^2\Bigr)^{1/2},\\
        \beta_F & :=\Bigl(\int_{\mathcal E}\sum_j\|F_j(e)\|_2^2\lambda(\mathrm de)\Bigr)^{1/2},
        & c_0 & :=\|C_0\|_2,\qquad
        e_0:=\Bigl(\int_{\mathcal E}\|E_0(e)\|_2^2\lambda(\mathrm de)\Bigr)^{1/2},\\
        a_1 & :=2\|S\|_2(\beta_B+c_0\beta_D+e_0\beta_F),
        & a_2 & :=\|S\|_2(\beta_D^2+\beta_F^2).
    \end{aligned}
\end{equation*}

\begin{lemma}
    \label{lem:certificate-stability}
    If \(a_1\delta+a_2\delta^2<\alpha\), then \(\bar B\subset\mathcal G\).
\end{lemma}
\begin{proof}
    Fix \(K\in\bar B\), and let \(\Delta K_j=K_j-K_j^{(0)}\), \(h=\|K-K^{(0)}\|_*\),
    \begin{equation*}
        \label{eq:cert-coefficient-perturbations}
        \begin{aligned}
            \Delta A & :=A_K-A_{K^{(0)}}=\sum_jB_j\Delta K_j,\\
            \Delta C & :=C_K-C_0=\sum_jD_j\Delta K_j,\\
            \Delta E(e) & :=E_K(e)-E_0(e)=\sum_jF_j(e)\Delta K_j.
        \end{aligned}
    \end{equation*}
    By the Cauchy--Schwarz inequality over the player blocks,
    \begin{equation*}
        \label{eq:cert-drift-diffusion-coefficient-bounds}
        \|\Delta A\|_2
        \le\sum_j\|B_j\|_2\|\Delta K_j\|_2\le \beta_Bh,
        \qquad
        \|\Delta C\|_2\le \beta_Dh.
    \end{equation*}
    \begin{equation*}
        \label{eq:cert-jump-coefficient-bound}
        \Bigl(\int_{\mathcal E}\|\Delta E(e)\|_2^2\lambda(\mathrm de)\Bigr)^{1/2}
        \le
        \Bigl(\int_{\mathcal E}\sum_j\|F_j(e)\|_2^2\lambda(\mathrm de)\Bigr)^{1/2}
        \Bigl(\sum_j\|\Delta K_j\|_2^2\Bigr)^{1/2}
        =\beta_Fh.
    \end{equation*}
    The drift part of \(\mathcal L_K(S)-\mathcal L_{K^{(0)}}(S)\) is \(\Delta A^\top S+S\Delta A\), so
    \begin{equation*}
        \label{eq:cert-drift-perturbation-bound}
        \|\Delta A^\top S+S\Delta A\|_2
        \le2\|S\|_2\beta_Bh.
    \end{equation*}
    For the diffusion term,
    \begin{equation*}
        \label{eq:cert-diffusion-expansion}
        C_K^\top SC_K-C_0^\top SC_0
        =C_0^\top S\Delta C+\Delta C^\top SC_0+\Delta C^\top S\Delta C,
    \end{equation*}
    and hence
    \begin{equation*}
        \label{eq:cert-diffusion-perturbation-bound}
        \|C_K^\top SC_K-C_0^\top SC_0\|_2
        \le2\|S\|_2c_0\beta_Dh+\|S\|_2\beta_D^2h^2.
    \end{equation*}
    Likewise, by the Cauchy--Schwarz inequality with respect to \(\lambda\),
    \begin{equation*}
        \label{eq:cert-jump-perturbation-bound}
        \begin{aligned}
            & \Bigl\|\int_{\mathcal E}
            \bigl[E_K(e)^\top SE_K(e)-E_0(e)^\top SE_0(e)\bigr]
            \lambda(\mathrm de)\Bigr\|_2\le
            2\|S\|_2e_0\beta_Fh+\|S\|_2\beta_F^2h^2.
        \end{aligned}
    \end{equation*}
    Combining the three bounds and the definitions of \(a_1,a_2\) gives
    \begin{equation}
        \label{eq:cert-lyapunov-perturbation}
        \|\mathcal L_K(S)-\mathcal L_{K^{(0)}}(S)\|_2\le a_1h+a_2h^2.
    \end{equation}
    Since the two matrices \(\mathcal L_K(S)\) and \(\mathcal L_{K^{(0)}}(S)\) in~\eqref{eq:cert-lyapunov-perturbation} are symmetric, their difference is bounded above by its spectral norm times \(I_n\). For \(h\le\delta\), combining~\eqref{eq:cert-center-lyapunov} and~\eqref{eq:cert-lyapunov-perturbation} gives
    \begin{equation*}
        \label{eq:cert-uniform-lyapunov-inequality}
        \mathcal L_K(S)
        \leq\mathcal L_{K^{(0)}}(S)
        +\|\mathcal L_K(S)-\mathcal L_{K^{(0)}}(S)\|_2 I_n
        \leq-[\alpha-a_1\delta-a_2\delta^2]I_n<0.
    \end{equation*}
    The strict Lyapunov inequality implies mean-square stability.
\end{proof}
\begin{remark}[Explicit admissible radius]
    \label{rem:cert-stability-radius}
    If \(a_2>0\), every radius satisfying
    \begin{equation*}
        \label{eq:cert-stability-radius}
        0<\delta<
        \frac{-a_1+\sqrt{a_1^2+4a_2\alpha}}{2a_2}
    \end{equation*}
    preserves mean-square stability. If \(a_2=0<a_1\), it is sufficient to require \(0<\delta<\alpha/a_1\); if \(a_1=a_2=0\), the stability condition holds for every \(\delta>0\).
\end{remark}

For \(K\in\mathcal G\), write \(\mathcal P(K)=(\mathcal P_1(K),\ldots,\mathcal P_N(K))\), and set \(\mathcal M_0:=\mathcal M(\mathcal P(K^{(0)}))\). Denote the smallest singular value of \(\mathcal M_0\) by \(\mu_0\), i.e., \(\mu_0:=\sigma_{\min}(\mathcal M_0)>0\).

On the stabilizing ball supplied by Lemma~\ref{lem:certificate-stability}, let \(L_{\mathcal M}>0\) be a bound satisfying
\begin{equation}
    \label{eq:cert-improvement-lipschitz-bound}
    \|\mathcal M(\mathcal P(K))-\mathcal M_0\|_2\le L_{\mathcal M}\|K-K^{(0)}\|_*,\qquad K\in\bar B.
\end{equation}
Such a bound can be estimated using \(D\mathcal P_i\) from~\eqref{eq:pi-evaluation-derivative}.

\begin{lemma}
    \label{lem:certificate-domain}
    Under the conditions of Lemma~\ref{lem:certificate-stability}, if \(L_{\mathcal M}\delta<\mu_0\), then \(\bar B\subset\mathcal D\) and
    \begin{equation*}
        \label{eq:cert-uniform-singular-value-bound}
        \sigma_{\min}\bigl(\mathcal M(\mathcal P(K))\bigr)\ge\mu_0-L_{\mathcal M}\delta>0,
        \qquad K\in\bar B.
    \end{equation*}
\end{lemma}
\begin{proof}
    By the singular-value perturbation inequality,
    \begin{equation*}
        \label{eq:cert-singular-value-perturbation}
        \sigma_{\min}\bigl(\mathcal M(\mathcal P(K))\bigr)
        \ge
        \sigma_{\min}(\mathcal M_0)
        -\|\mathcal M(\mathcal P(K))-\mathcal M_0\|_{2}.
    \end{equation*}
    Using \eqref{eq:cert-improvement-lipschitz-bound} and \(\|K-K^{(0)}\|_{*}\le \delta\) gives
    \begin{equation*}
        \label{eq:cert-singular-value-lower-bound}
        \sigma_{\min}\bigl(\mathcal M(\mathcal P(K))\bigr)
        \ge
        \mu_0-L_{\mathcal M}\delta>0.
    \end{equation*}
    Hence \(\mathcal M(\mathcal P(K))\) is nonsingular throughout the ball.
\end{proof}

For a candidate derivative bound \(L_{\mathcal T}>0\), define
\begin{equation}
    \label{eq:cert-contraction-constants}
    q_0:=\|D\mathcal T(K^{(0)})\|_*,\qquad
    q:=q_0+L_{\mathcal T}\delta,\qquad
    \eta:=\|\mathcal T(K^{(0)})-K^{(0)}\|_*.
\end{equation}
The following assumption collects the local stability, domain, and derivative conditions needed for a verifiable convergence certificate.
\begin{assumption}
    \label{ass:local-certificate}
    For the selected \(K^{(0)}\in\mathcal D\) and radius \(\delta>0\), suppose there exist \(S\in\mathcal S_{++}^{n}\), \(\alpha>0\), and constants \(L_{\mathcal M},L_{\mathcal T}>0\) such that:
    \begin{enumerate}
        \item[(a)] \eqref{eq:cert-center-lyapunov} holds, and the stability margin satisfies
            \begin{equation*}
                \label{eq:cert-stability-condition}
                a_1\delta+a_2\delta^2<\alpha;
            \end{equation*}
        \item[(b)] the nonsingularity margin satisfies
            \begin{equation*}
                \label{eq:cert-nonsingularity-condition}
                L_{\mathcal M}\delta<\mu_0;
            \end{equation*}
        \item[(c)] \(\mathcal T\) is continuously differentiable on \(\bar B\), and the derivative variation satisfies, for every \(K\in\bar B\),
            \begin{equation*}
                \label{eq:cert-derivative-lipschitz-bound}
                \|D\mathcal T(K)-D\mathcal T(K^{(0)})\|_*
                \le L_{\mathcal T}\|K-K^{(0)}\|_*.
            \end{equation*}
            Moreover, \(q<1\);
        \item[(d)] the initial residual satisfies
            \begin{equation*}
                \label{eq:cert-residual-condition}
                \eta\le(1-q)\delta.
            \end{equation*}
    \end{enumerate}
\end{assumption}

We now state the local convergence certificate. It uses these checkable bounds, which are all centered at \(K^{(0)}\), to establish invariance and contraction for the policy-iteration map.

\begin{theorem}[Verifiable local convergence certificate]
    \label{thm:local-certificate}
    Suppose Assumption~\ref{ass:local-certificate} holds. Then:
    \begin{enumerate}
        \item[(i)] \(\bar B\subset\mathcal D\) and \(\mathcal T(\bar B)\subset\bar B\);
        \item[(ii)] \(\mathcal T\) is a \(q\)-contraction on \(\bar B\) and has a unique fixed point \(K^*\in\bar B\);
        \item[(iii)] Algorithm~\ref{alg:policy-iteration} is well defined, every iterate is mean-square stabilizing, and
            \begin{equation*}
                \label{eq:cert-geometric-convergence-bound}
                \|K^{(l)}-K^*\|_*
                \le\frac{q^l\eta}{1-q},\qquad l\ge0,
            \end{equation*}
            where \(K^*\) is a closed-loop Nash equilibrium.
    \end{enumerate}
\end{theorem}

\begin{proof}
    \begin{enumerate}
        \item[(i)] By Assumption~\ref{ass:local-certificate}(a) and Lemma~\ref{lem:certificate-stability}, \(\bar B\subset\mathcal G\). Assumption~\ref{ass:local-certificate}(b) and Lemma~\ref{lem:certificate-domain} ensure that the policy-improvement matrix is nonsingular throughout \(\bar B\), so \(\bar B\subset\mathcal D\). For \(K\in\bar B\), Assumption~\ref{ass:local-certificate}(c) and~\eqref{eq:cert-contraction-constants} give
            \begin{equation*}
                \label{eq:cert-uniform-contraction-bound}
                \|D\mathcal T(K)\|_*
                \le\|D\mathcal T(K^{(0)})\|_*+L_{\mathcal T}\|K-K^{(0)}\|_*
                \le q<1.
            \end{equation*}
            Hence the mean-value inequality yields \(\|\mathcal T(K)-\mathcal T(K')\|_*\le q\|K-K'\|_*\) for all \(K,K'\in\bar B\). In particular, Assumption~\ref{ass:local-certificate}(d) implies
            \begin{equation*}
                \label{eq:cert-ball-invariance-bound}
                \|\mathcal T(K)-K^{(0)}\|_*
                \le\|\mathcal T(K)-\mathcal T(K^{(0)})\|_*+\eta
                \le q\delta+\eta\le\delta,
            \end{equation*}
            and therefore \(\mathcal T(\bar B)\subset\bar B\).

        \item[(ii)] The ball \(\bar B\) is complete, and part~(i) shows that \(\mathcal T\) maps it into itself and is a contraction there. Banach's fixed-point theorem gives a unique fixed point \(K^*\in\bar B\).

        \item[(iii)] Since \(K^{(0)}\in\bar B\) and \(\mathcal T(\bar B)\subset\bar B\subset\mathcal D\), every policy-iteration step is well defined and mean-square stabilizing. Moreover,
            \begin{equation*}
                \label{eq:cert-successive-increment-bound}
                \|K^{(l+1)}-K^{(l)}\|_*
                \le q^l\|K^{(1)}-K^{(0)}\|_*=q^l\eta.
            \end{equation*}
            Summing the geometric tail gives \(\|K^{(l)}-K^*\|_*\le q^l\eta/(1-q)\) for all \(l\ge0\).
    \end{enumerate}
\end{proof}

\begin{remark}
	If Assumption~\ref{ass:local-certificate} fails at \(K^{(0)}\) but holds on a ball centered at a later iterate \(K^{(l)}\), then Theorem~\ref{thm:local-certificate} can be applied with \(K^{(l)}\)	as the initial feedback profile.
\end{remark}

\subsection{Totally Model-free Off-policy Q-learning algorithm}
\label{subsec:off-policy-q-learning}
In this subsection, we express policy evaluation and simultaneous policy improvement through each player's Q-function. Behavior-policy data yield an off-policy regression that can be reused across target-policy updates without system coefficients in the learning step. We show that the exact population regression recovers the same value matrices and updated gains as policy iteration under quantitative persistence of excitation.

For each player \(i \in \mathbb N\), let \(Y = [x^{\top},u_{1}^{\top},\ldots,u_{N}^{\top}]^{\top} \in \mathbb{R}^{d},\ d:= n+\sum_{j = 1}^{N}m_{j}\), and define the Q-function associated with \(P_{i} \in \mathcal{S}_{++}^{n}\) by
\begin{equation}
    \label{eq:q-function-definition}
    \mathcal{Q}_{i}(x,u_{1},\ldots,u_{N};P_{i})
    := Y^{\top}\mathbf{Q}_{i}(P_{i})Y
    = x^{\top} P_{i}x+\mathcal{H}_{i}(x,u_{1},\ldots,u_{N};P_{i}).
\end{equation}
The matrix \(\mathbf{Q}_{i}(P_{i}) \in \mathcal{S}^{d}\) is partitioned as
\begin{equation}
    \label{eq:q-matrix-partition}
    \mathbf{Q}_{i}(P_{i}) =
    \begin{bmatrix}
        Q_{i,xx}     & Q_{i,xu_{1}}     & \cdots & Q_{i,xu_{N}}     \\
        Q_{i,u_{1}x} & Q_{i,u_{1}u_{1}} & \cdots & Q_{i,u_{1}u_{N}} \\
        \vdots       & \vdots           & \ddots & \vdots           \\
        Q_{i,u_{N}x} & Q_{i,u_{N}u_{1}} & \cdots & Q_{i,u_{N}u_{N}}
    \end{bmatrix},
\end{equation}
where
\begin{equation}
    \label{eq:q-matrix-blocks}
    \begin{aligned}
        Q_{i,xx}
        & =P_{i}+A^{\top} P_{i}+P_{i}A+C^{\top} P_{i}C+N_{i}
        +\int_{\mathcal{E}}E(e)^{\top} P_{i}E(e)\lambda(\mathrm{d}e),\\
        Q_{i,xu_{j}}
        & =P_{i}B_{j}+C^{\top} P_{i}D_{j}
        +\int_{\mathcal{E}}E(e)^{\top} P_{i}F_{j}(e)\lambda(\mathrm{d}e),
        \qquad Q_{i,u_{j}x} = Q_{i,xu_{j}}^{\top},\\
        Q_{i,u_{j}u_{j}}
        & =R_{ij}
        +D_{j}^{\top} P_{i}D_{j}
        +\int_{\mathcal{E}}F_{j}(e)^{\top} P_{i}F_{j}(e)\lambda(\mathrm{d}e),\\
        Q_{i,u_{j}u_{k}}
        & =D_{j}^{\top} P_{i}D_{k}
        +\int_{\mathcal{E}}F_{j}(e)^{\top} P_{i}F_{k}(e)\lambda(\mathrm{d}e),
        \qquad j,k \in \mathbb N,\quad j \neq k.
    \end{aligned}
\end{equation}
Using~\eqref{eq:q-matrix-partition} and~\eqref{eq:q-matrix-blocks}, combined with the definition of the Hamiltonian function, we can rewrite the Q-function~\eqref{eq:q-function-definition} as
\begin{equation}
    \label{eq:q-function-expansion}
    \begin{aligned}
        \mathcal{Q}_{i}(x,u_{1},\ldots,u_{N};P_{i})
        & =x^{\top}\left[P_{i}+A^{\top} P_{i}+P_{i}A+C^{\top} P_{i}C+N_{i}
            +\int_{\mathcal{E}}E(e)^{\top} P_{i}E(e)\lambda(\mathrm{d}e)\right]x\\
        & \quad+2\sum_{j = 1}^{N}x^{\top}\left[P_{i}B_{j}+C^{\top} P_{i}D_{j}
            +\int_{\mathcal{E}}E(e)^{\top} P_{i}F_{j}(e)\lambda(\mathrm{d}e)\right]u_{j}\\
        & \quad+\sum_{j = 1}^{N}u_{j}^{\top}\left[R_{ij}+D_{j}^{\top} P_{i}D_{j}
            +\int_{\mathcal{E}}F_{j}(e)^{\top} P_{i}F_{j}(e)\lambda(\mathrm{d}e)\right]u_{j}\\
        & \quad+2\sum_{1 \leq j < k \leq N}u_{j}^{\top}\left[D_{j}^{\top} P_{i}D_{k}
            +\int_{\mathcal{E}}F_{j}(e)^{\top} P_{i}F_{k}(e)\lambda(\mathrm{d}e)\right]u_{k}\\
        & =x^{\top} P_{i}x+\mathcal{H}_{i}(x,u_{1},\ldots,u_{N};P_{i}).
    \end{aligned}
\end{equation}
At iteration \(l\), let \(P_{i}^{(l+1)}\) be the solution of the Lyapunov recursion~\eqref{eq:pi-lyapunov-recursion} under \((K_{1}^{(l)}, \dots, K_{N}^{(l)})\). For simplicity, denote
\begin{equation*}
    \label{eq:q-iteration-notation}
    \mathcal{Q}_{i}^{(l+1)}(x,u_{1},\ldots,u_{N})
    := \mathcal{Q}_{i}(x,u_{1},\ldots,u_{N};P_{i}^{(l+1)})
    = Y^{\top}\mathbf{Q}_{i}^{(l+1)}Y,
\end{equation*}
Substituting \(u_{j}^{(l)} = K_{j}^{(l)}x\), \(j \in \mathbb N\), we obtain
\begin{equation}
    \label{eq:q-target-value}
    \mathcal{Q}_{i}^{(l+1)}
    \bigl(x,K_{1}^{(l)}x,\ldots,K_{N}^{(l)}x\bigr)
    = x^{\top} P_{i}^{(l+1)}x
    =:V_{i}^{(l+1)}(x),\qquad i \in \mathbb N.
\end{equation}

We now introduce the behavior policy \(\widetilde{u}_{j}(\cdot)\) used to generate the data. The state trajectory generated by the behavior policy satisfies
\begin{equation}
    \label{eq:q-behavior-dynamics}
    \begin{aligned}
        \mathrm{d}X(s)
        & =\left[A_{l}X(s)
            +\sum_{j = 1}^{N}B_{j}\bigl(\widetilde{u}_{j}(s)-K_{j}^{(l)}X(s)\bigr)\right]\mathrm{d}s\\
        & \quad+\left[C_{l}X(s)
            +\sum_{j = 1}^{N}D_{j}\bigl(\widetilde{u}_{j}(s)-K_{j}^{(l)}X(s)\bigr)\right]\mathrm{d}W(s)\\
        & \quad+\int_{\mathcal{E}}\left[E_{l}(e)X(s-)
            +\sum_{j = 1}^{N}F_{j}(e)\bigl(\widetilde{u}_{j}(s)-K_{j}^{(l)}X(s-)\bigr)\right]
        \widetilde{N}(\mathrm{d}s,\mathrm{d}e).
    \end{aligned}
\end{equation}
Let \(Y_{l}(s):= [X(s)^{\top},(K_{1}^{(l)}X(s))^{\top},\ldots,(K_{N}^{(l)}X(s))^{\top}]^{\top}, \widetilde{Y}(s):= [X(s)^{\top},\widetilde{u}_{1}(s)^{\top},\ldots,\widetilde{u}_{N}(s)^{\top}]^{\top}\). By equation~\eqref{eq:q-function-definition}, the target-policy and behavior-policy evaluations are, respectively,
\begin{equation}
    \label{eq:q-target-evaluation}
    \mathcal{Q}_{i}^{(l+1)}
    \bigl(X(s),K_{1}^{(l)}X(s),\ldots,K_{N}^{(l)}X(s)\bigr)
    = Y_{l}(s)^{\top}\mathbf{Q}_{i}^{(l+1)}Y_{l}(s),
\end{equation}
and
\begin{equation}
    \label{eq:q-behavior-evaluation}
    \mathcal{Q}_{i}^{(l+1)}
    \bigl(X(s),\widetilde{u}_{1}(s),\ldots,\widetilde{u}_{N}(s)\bigr)
    = \widetilde{Y}(s)^{\top}\mathbf{Q}_{i}^{(l+1)}\widetilde{Y}(s).
\end{equation}
Then, for every \(t \geq 0\), \(\Delta t > 0\), and \(i \in \mathbb N\), \(X(\cdot)\) is generated by the fixed behavior SDE~\eqref{eq:q-behavior-dynamics}. We have the following equation
\begin{equation}
    \label{eq:q-integral-identity}
    \begin{aligned}
        & \mathcal{Q}_{i}^{(l+1)}
        \bigl(X(t),K_{1}^{(l)}X(t),\ldots,K_{N}^{(l)}X(t)\bigr)\\
        & \quad-\mathbb{E}^{\mathcal{F}_{t}}\!\left[
            \mathcal{Q}_{i}^{(l+1)}
            \bigl(X(t+\Delta t),K_{1}^{(l)}X(t+\Delta t),\ldots,
            K_{N}^{(l)}X(t+\Delta t)\bigr)\right]\\
        & \quad+\mathbb{E}^{\mathcal{F}_{t}}\!\int_{t}^{t+\Delta t}
        \left[
            \mathcal{Q}_{i}^{(l+1)}\bigl(X(s),\widetilde{u}_{1}(s),\ldots,\widetilde{u}_{N}(s)\bigr)
            -\mathcal{Q}_{i}^{(l+1)}
            \bigl(X(s),K_{1}^{(l)}X(s),\ldots,K_{N}^{(l)}X(s)\bigr)
            \right]\mathrm{d}s\\
        & =\mathbb{E}^{\mathcal{F}_{t}}\!\int_{t}^{t+\Delta t}
        \left[X(s)^{\top} N_{i}X(s)+
            \sum_{j = 1}^{N}\widetilde{u}_{j}(s)^{\top} R_{ij}\widetilde{u}_{j}(s)\right]\mathrm{d}s.
    \end{aligned}
\end{equation}
Substituting~\eqref{eq:q-target-evaluation} and~\eqref{eq:q-behavior-evaluation} into~\eqref{eq:q-integral-identity}, we obtain
\begin{equation}
    \label{eq:q-regression-identity}
    \begin{aligned}
        & Y_{l}(t)^{\top}\mathbf{Q}_{i}^{(l+1)}Y_{l}(t)
        -\mathbb{E}^{\mathcal{F}_{t}}\!\left[
            Y_{l}(t+\Delta t)^{\top}\mathbf{Q}_{i}^{(l+1)}Y_{l}(t+\Delta t)\right]\\
        & \quad+\mathbb{E}^{\mathcal{F}_{t}}\!\int_{t}^{t+\Delta t}
        \left[\widetilde{Y}(s)^{\top}\mathbf{Q}_{i}^{(l+1)}\widetilde{Y}(s)
            -Y_{l}(s)^{\top}\mathbf{Q}_{i}^{(l+1)}Y_{l}(s)\right]\mathrm{d}s\\
        & =\mathbb{E}^{\mathcal{F}_{t}}\!\int_{t}^{t+\Delta t}
        \left[X(s)^{\top} N_{i}X(s)
            +\sum_{j = 1}^{N}\widetilde{u}_{j}(s)^{\top} R_{ij}\widetilde{u}_{j}(s)\right]\mathrm{d}s,
        \qquad i \in \mathbb N.
    \end{aligned}
\end{equation}
Once \(\mathbf{Q}_{i}^{(l+1)}\) is identified, set \(\mathbf{K}^{(l)}:= [I_{n}^{\top},K_{1}^{(l)\top},\ldots,K_{N}^{(l)\top}]^{\top} \in \mathbb{R}^{d\times n}\). The value matrices are recovered through
\begin{equation}
    \label{eq:q-value-recovery}
    P_{i}^{(l+1)}
    = \mathbf{K}^{(l)\top}\mathbf{Q}_{i}^{(l+1)}\mathbf{K}^{(l)},
    \qquad i \in \mathbb N,
\end{equation}
and all feedback gains are updated by
\begin{equation}
    \label{eq:q-gain-update}
    \begin{bmatrix}K_{1}^{(l+1)}\\\vdots\\K_{N}^{(l+1)}\end{bmatrix}
    = -\begin{bmatrix}
        Q_{1,u_{1}u_{1}}^{(l+1)} & \cdots & Q_{1,u_{1}u_{N}}^{(l+1)} \\
        \vdots                   & \ddots & \vdots                   \\
        Q_{N,u_{N}u_{1}}^{(l+1)} & \cdots & Q_{N,u_{N}u_{N}}^{(l+1)}
    \end{bmatrix}^{-1}
    \begin{bmatrix}Q_{1,u_{1}x}^{(l+1)}\\\vdots\\Q_{N,u_{N}x}^{(l+1)}\end{bmatrix}.
\end{equation}

Using the Q-function representation, we introduce a totally model-free off-policy Q-learning algorithm that reuses the same behavior-policy data across iterations.

\begin{algorithm}[H]
    \caption{Totally model-free off-policy Q-learning for the \(N\)-player Problem (LQ NZSSDG)}
    \label{alg:off-policy-q-learning}
    \textbf{Initialize:} Select an admissible initial target feedback profile \(K^{(0)} = (K_{1}^{(0)},\ldots,K_{N}^{(0)}) \in \mathcal{D}\).
    \\[2mm]
    \textbf{Data Collection:} Apply behavior policies \(\widetilde{u}_{j}\), \(j \in \mathbb N\), to the unknown stochastic system and collect one batch \(\{X(s),\widetilde{u}_{1}(s),\ldots,\widetilde{u}_{N}(s)\}_{s \in [t_{0},t_{g}]}\).
    \\[2mm]
    \textbf{For \(l = 0,1,2,\ldots\):}
    \begin{itemize}
        \item \textbf{Policy Evaluation:} Solve for \(\mathbf{Q}_{i}^{(l+1)}\) using~\eqref{eq:q-regression-identity}, and recover \(P_{i}^{(l+1)}\) by~\eqref{eq:q-value-recovery}.
        \item \textbf{Policy Improvement:} Update all gains by~\eqref{eq:q-gain-update}.
    \end{itemize}
    \textbf{Stop when} \(\sum_{i = 1}^{N}\|P_{i}^{(l+1)}-P_{i}^{(l)}\|_F < \varepsilon\) for \(l \geq 1\).
\end{algorithm}
\begin{lemma}
    \label{lem:q-policy-evaluation}
    For every \(K^{(l)} \in \mathcal{D}\), the matrices \(\mathbf{Q}_{i}^{(l+1)} = \mathbf{Q}_{i}(P_{i}^{(l+1)})\) constructed from the solution of~\eqref{eq:pi-lyapunov-recursion} satisfy~\eqref{eq:q-regression-identity} in the policy evaluation step of Algorithm~\ref{alg:off-policy-q-learning} for all \(i \in \mathbb N\).
\end{lemma}

\begin{proof}
    Applying Itô's formula to \(X(s)^{\top} P_{i}^{(l+1)}X(s)\) along the behavior dynamics~\eqref{eq:q-behavior-dynamics} gives
    \begin{equation}
        \label{eq:q-behavior-ito}
        \begin{aligned}
            & \mathrm{d}\bigl(X(s)^{\top} P_{i}^{(l+1)}X(s)\bigr)\\
            & =\Bigl\{X(s)^{\top} P_{i}^{(l+1)}\Bigl(AX(s)+\sum_{j = 1}^{N}B_{j}\widetilde{u}_{j}(s)\Bigr)+\Bigl(AX(s)+\sum_{j = 1}^{N}B_{j}\widetilde{u}_{j}(s)\Bigr)^{\top} P_{i}^{(l+1)}X(s)\\
            & \quad+\Bigl(CX(s)+\sum_{j = 1}^{N}D_{j}\widetilde{u}_{j}(s)\Bigr)^{\top} P_{i}^{(l+1)}\Bigl(CX(s)+\sum_{j = 1}^{N}D_{j}\widetilde{u}_{j}(s)\Bigr)\\
            & \quad+\int_{\mathcal{E}}\Bigl(E(e)X(s)+\sum_{j = 1}^{N}F_{j}(e)\widetilde{u}_{j}(s)\Bigr)^{\top} P_{i}^{(l+1)}\Bigl(E(e)X(s)+\sum_{j = 1}^{N}F_{j}(e)\widetilde{u}_{j}(s)\Bigr)\lambda(\mathrm{d}e)\Bigr\}\mathrm{d}s\\
            & \quad+\{\cdots\}\,\mathrm{d}W(s) + \{\cdots\}\,\widetilde{N}(\mathrm{d}s,\mathrm{d}e).
        \end{aligned}
    \end{equation}
    By~\eqref{eq:q-function-expansion}, we can obtain the Q-function under the behavior policy
    \begin{equation}
        \label{eq:q-behavior-expansion}
        \begin{aligned}
            & \mathcal{Q}_{i}^{(l+1)}\bigl(X(s),\widetilde{u}_{1}(s),\ldots,\widetilde{u}_{N}(s)\bigr)\\
            & \quad = X(s)^{\top} P_{i}^{(l+1)} X(s) + \mathcal{H}_{i}(X(s),\widetilde{u}_{1}(s),\ldots,\widetilde{u}_{N}(s);P_{i}^{(l+1)})\\
            & \quad = X(s)^{\top} P_{i}^{(l+1)} X(s) + X(s)^{\top} P_{i}^{(l+1)}\Bigl(AX(s)+\sum_{j = 1}^{N}B_{j}\widetilde{u}_{j}(s)\Bigr)+ \Bigl(AX(s)+\sum_{j = 1}^{N}B_{j}\widetilde{u}_{j}(s)\Bigr)^{\top} P_{i}^{(l+1)}X(s)\\
            & \qquad+\Bigl(CX(s)+\sum_{j = 1}^{N}D_{j}\widetilde{u}_{j}(s)\Bigr)^{\top} P_{i}^{(l+1)}\Bigl(CX(s)+\sum_{j = 1}^{N}D_{j}\widetilde{u}_{j}(s)\Bigr)\\
            & \qquad+\int_{\mathcal{E}}\Bigl(E(e)X(s)+\sum_{j = 1}^{N}F_{j}(e)\widetilde{u}_{j}(s)\Bigr)^{\top} P_{i}^{(l+1)}\Bigl(E(e)X(s)+\sum_{j = 1}^{N}F_{j}(e)\widetilde{u}_{j}(s)\Bigr)\lambda(\mathrm{d}e)\\
            & \qquad+X(s)^{\top} N_{i}X(s)+\sum_{j = 1}^{N}\widetilde{u}_{j}(s)^{\top} R_{ij}\widetilde{u}_{j}(s).
        \end{aligned}
    \end{equation}
    Combining~\eqref{eq:q-behavior-ito} and~\eqref{eq:q-behavior-expansion} yields
    \begin{equation}
        \label{eq:q-value-differential}
        \begin{aligned}
            \mathrm{d}\bigl(X(s)^{\top} P_{i}^{(l+1)}X(s)\bigr)
            & =\Bigl[\mathcal{Q}_{i}^{(l+1)}
                \bigl(X(s),\widetilde{u}_{1}(s),\ldots,\widetilde{u}_{N}(s)\bigr)
                -X(s)^{\top} P_{i}^{(l+1)}X(s)\\
                & \qquad-X(s)^{\top} N_{i}X(s)
                -\sum_{j = 1}^{N}\widetilde{u}_{j}(s)^{\top} R_{ij}\widetilde{u}_{j}(s)
                \Bigr]\mathrm{d}s\\
            & \qquad+\{\cdots\}\,\mathrm{d}W(s) + \{\cdots\}\,\widetilde{N}(\mathrm{d}s,\mathrm{d}e),
        \end{aligned}
    \end{equation}
    Integrating both sides of~\eqref{eq:q-value-differential} over \([t,t+\Delta t]\) and taking the conditional expectation \(\mathbb{E}^{\mathcal{F}_{t}}\), we obtain
    \begin{equation}
        \label{eq:q-integrated-value-identity}
        \begin{aligned}
            \mathbb{E}^{\mathcal{F}_{t}}\! & \left[
                X(t+\Delta t)^{\top} P_{i}^{(l+1)}X(t+\Delta t)\right]-X(t)^{\top} P_{i}^{(l+1)}X(t)\\
            & =\mathbb{E}^{\mathcal{F}_{t}}\!\int_{t}^{t+\Delta t}
            \Bigl[\mathcal{Q}_{i}^{(l+1)}
                \bigl(X(s),\widetilde{u}_{1}(s),\ldots,\widetilde{u}_{N}(s)\bigr)
                -X(s)^{\top} P_{i}^{(l+1)}X(s)\\
                & \hspace{40mm}-X(s)^{\top} N_{i}X(s)
                -\sum_{j = 1}^{N}\widetilde{u}_{j}(s)^{\top} R_{ij}\widetilde{u}_{j}(s)
                \Bigr]\mathrm{d}s.
        \end{aligned}
    \end{equation}
    Finally, substituting~\eqref{eq:q-target-value} into~\eqref{eq:q-integrated-value-identity} and rearranging the terms yields~\eqref{eq:q-integral-identity}.
\end{proof}
\begin{remark}
    \label{rem:q-identifiability}
    Lemma~\ref{lem:q-policy-evaluation} establishes the consistency between the Lyapunov recursion and the policy evaluation step in Algorithm~\ref{alg:off-policy-q-learning}. Furthermore, the solution \(P_{i}^{(l+1)}\) of Algorithm~\ref{alg:policy-iteration} is a solution to~\eqref{eq:q-integral-identity}. It should be noted that this lemma does not guarantee the uniqueness of the solution. To ensure that~\eqref{eq:q-integral-identity} has a unique solution, the persistence of excitation (PE) condition must be satisfied. Only under this condition can Algorithm~\ref{alg:off-policy-q-learning} uniquely recover the same \(P_{i}^{(l+1)}\) and \(K_{i}^{(l+1)}\).
\end{remark}
Next, we introduce the vectorization and Kronecker-product notation used below. For \(A = (a_{ij}) \in \mathcal{S}^{m}\), define \(\operatorname{svec}(A):= [a_{11},2a_{12},\ldots,2a_{1m},a_{22},2a_{23},\ldots,2a_{m-1,m},a_{mm}]^{\top} \in \mathbb{R}^{m(m+1)/2}\). For \(y \in \mathbb{R}^{d}\), we have \(y^{\top} \mathbf{Q} y = \bar{y}^{\top}\operatorname{svec}(\mathbf{Q})\) for every \(\mathbf{Q} \in \mathcal{S}^{d}\). Let \(D_{d} \in \mathbb{R}^{d^{2}\times d(d+1)/2}\) be the matrix satisfying \(y\otimes y = D_{d}\bar{y}\). Let \(A\otimes B\) denote the Kronecker product of matrices \(A\) and \(B\). Then \(\operatorname{vec}(ABC) = (C^{\top}\otimes A)\operatorname{vec}(B)\) and \((A\otimes B)(C\otimes D) = (AC)\otimes(BD)\).

Let \(0 \leq t_{0} < t_{1} < \cdots < t_{g}\), where \(g\) is a positive integer. For \(j \in \mathbb N\), we define matrices \(\eta_{y_{l}y_{l}},\mathcal{I}_{\widetilde{y} \widetilde{y}},\mathcal{I}_{y_{l}y_{l}} \in \mathbb{R}^{g\times\frac{d(d+1)}{2}},\ \mathcal{I}_{xx} \in \mathbb{R}^{g\times n^{2}},\ \mathcal{I}_{\widetilde{u}_{j}\widetilde{u}_{j}} \in \mathbb{R}^{g\times m_{j}^{2}}\) as follows
\begin{equation*}
    \label{eq:q-data-matrices}
    \begin{aligned}
        \eta_{y_{l}y_{l}}
        & :=\left[\overline{Y_{l}(t_{0})}^{\top}-\mathbb{E}^{\mathcal{F}_{t_{0}}}\left[\overline{Y_{l}(t_{1})}^{\top}\right],\ldots,\overline{Y_{l}(t_{g-1})}^{\top}-\mathbb{E}^{\mathcal{F}_{t_{g-1}}}\left[\overline{Y_{l}(t_{g})}^{\top}\right]\right]^{\top},\\[1mm]
        \mathcal{I}_{\widetilde{y} \widetilde{y}}
        & :=
        \left[
            \mathbb{E}^{\mathcal{F}_{t_{0}}}
            \int_{t_{0}}^{t_{1}}
            \overline{\widetilde{Y}(s)}^{\top}\,\mathrm{d}s,
            \ldots,
            \mathbb{E}^{\mathcal{F}_{t_{g-1}}}
            \int_{t_{g-1}}^{t_{g}}
            \overline{\widetilde{Y}(s)}^{\top}\,\mathrm{d}s
            \right]^{\top},\\[1mm]
        \mathcal{I}_{y_{l}y_{l}}
        & :=
        \left[
            \mathbb{E}^{\mathcal{F}_{t_{0}}}
            \int_{t_{0}}^{t_{1}}
            \overline{Y_{l}(s)}^{\top}\,\mathrm{d}s,
            \ldots,
            \mathbb{E}^{\mathcal{F}_{t_{g-1}}}
            \int_{t_{g-1}}^{t_{g}}
            \overline{Y_{l}(s)}^{\top}\,\mathrm{d}s
            \right]^{\top},\\[1mm]
        \mathcal{I}_{\widetilde{u}_{j}\widetilde{u}_{j}}
        & :=
        \left[
            \mathbb{E}^{\mathcal{F}_{t_{0}}}
            \int_{t_{0}}^{t_{1}}
            \bigl(\widetilde{u}_{j}(s)^{\top}\otimes \widetilde{u}_{j}(s)^{\top}\bigr)\,\mathrm{d}s,
            \ldots,
            \mathbb{E}^{\mathcal{F}_{t_{g-1}}}
            \int_{t_{g-1}}^{t_{g}}
            \bigl(\widetilde{u}_{j}(s)^{\top}\otimes \widetilde{u}_{j}(s)^{\top}\bigr)\,\mathrm{d}s
            \right]^{\top},\\[1mm]
        \mathcal{I}_{xx}
        & :=
        \left[
            \mathbb{E}^{\mathcal{F}_{t_{0}}}
            \int_{t_{0}}^{t_{1}}
            \bigl(X(s)^{\top}\otimes X(s)^{\top}\bigr)\,\mathrm{d}s,
            \ldots,
            \mathbb{E}^{\mathcal{F}_{t_{g-1}}}
            \int_{t_{g-1}}^{t_{g}}
            \bigl(X(s)^{\top}\otimes X(s)^{\top}\bigr)\,\mathrm{d}s
            \right]^{\top}.
    \end{aligned}
\end{equation*}
Note that \(g\) should be no less than the number of unknown parameters of \(\mathbf{Q}_{i}^{(l+1)}\), i.e., \(g \geq \frac{d(d+1)}{2}\). With these symbols,~\eqref{eq:q-regression-identity} can be rewritten as
\begin{equation}
    \label{eq:q-linear-regression}
    \Phi^{(l)}\theta_{i}^{(l+1)} = \Pi_{i}, \qquad i \in \mathbb N,
\end{equation}
where \(\theta_{i}^{(l+1)} = \operatorname{svec}(\mathbf{Q}_{i}^{(l+1)})\), \(\Phi^{(l)} \in \mathbb{R}^{g\times\frac{d(d+1)}{2}}\) and \(\Pi_{i} \in \mathbb{R}^{g}\) are defined as
\begin{equation*}
    \label{eq:q-design-matrix}
    \Phi^{(l)}:= \eta_{y_{l}y_{l}}+\mathcal{I}_{\widetilde{y} \widetilde{y}}-\mathcal{I}_{y_{l}y_{l}},
\end{equation*}
\begin{equation*}
    \label{eq:q-regression-response}
    \Pi_{i}:= \mathcal{I}_{xx}\operatorname{vec}(N_{i})+\sum_{j = 1}^{N}\mathcal{I}_{\widetilde{u}_{j}\widetilde{u}_{j}}\operatorname{vec}(R_{ij}).
\end{equation*}
Notice that if \(\Phi^{(l)}\) has full column rank,~\eqref{eq:q-linear-regression} can be directly solved as follows:
\begin{equation*}
    \label{eq:q-least-squares-solution}
    \theta_{i}^{(l+1)} = \left(\Phi^{(l)\top} \Phi^{(l)}\right)^{-1} \Phi^{(l)\top} \Pi_{i}, \qquad i \in \mathbb N.
\end{equation*}

\begin{assumption}[Quantitative persistence of excitation]
    \label{ass:quantitative-excitation}
    Fix \(g\ge d(d+1)/2\). For every \(K^{(l)}\in\bar B\), \(l\ge0\), there exists a constant \(\kappa>0\) such that
    \begin{equation*}
        \label{eq:q-excitation-bound}
        \sigma_{\min}\!\left(\Phi^{(l)}\right)\ge\kappa.
    \end{equation*}
\end{assumption}

\begin{remark}
    \label{rem:q-excitation-design}
    Since \(\kappa>0\), Assumption~\ref{ass:quantitative-excitation} implies the full-column-rank condition
    \begin{equation}
        \label{eq:q-full-column-rank}
        \operatorname{rank}\Phi^{(l)}
        =\frac{d(d+1)}{2},
        \qquad l\ge0.
    \end{equation}
    It is stronger than requiring full column rank at each iterate, because it imposes a uniform positive lower bound over all admissible target profiles in \(\bar B\). The behavior policies are designed as \(\widetilde u_j(s)=K_j^{(0)}X(s)+e_j(s)\), \(j\in\mathbb N\), where \(e_j\) are bounded probing noises. It should be emphasized that \(e_{j}\) is not merely introduced as a general exploration noise. General exploration noise mainly perturbs the control input to enrich the trajectories, whereas the probing signals here are imposed to guarantee that \(\Phi^{(l)}\) has full column rank. Since the target policies satisfy \(u_{j}^{(l)} = K_{j}^{(l)}X\), the state and input data may become linearly dependent. The additional signals \(e_{j}\) break this dependence and ensure that the Q-function matrix can be uniquely identified from off-policy data.
\end{remark}

\begin{lemma}
    \label{lem:q-population-equivalence}
    Fix \(l\ge0\) and \(K^{(l)}\in\bar B\). Under Assumption~\ref{ass:quantitative-excitation}, Algorithm~\ref{alg:off-policy-q-learning} and Algorithm~\ref{alg:policy-iteration} produce the same \(P_i^{(l+1)}\) and the same updated feedback profile \(K^{(l+1)}\).
\end{lemma}

\begin{proof}
    Fix the common current feedback profile \((K_{1}^{(l)},\ldots,K_{N}^{(l)})\), and let \(P_{i}^{(l+1)}\) be the solution of the Lyapunov recursion~\eqref{eq:pi-lyapunov-recursion}. Therefore, we can construct \(\mathbf{Q}_{i}^{(l+1)} = \mathbf{Q}_{i}(P_{i}^{(l+1)})\) by~\eqref{eq:q-matrix-blocks}. Lemma~\ref{lem:q-policy-evaluation} indicates that \(\mathbf{Q}_{i}^{(l+1)}\) can be solved by~\eqref{eq:q-linear-regression}. Assumption~\ref{ass:quantitative-excitation} implies~\eqref{eq:q-full-column-rank}, so this solution is unique. Algorithm~\ref{alg:off-policy-q-learning} therefore identifies the same \(\mathbf{Q}_{i}^{(l+1)}\) and recovers the same value matrix \(P_{i}^{(l+1)}\) produced by Algorithm~\ref{alg:policy-iteration}.

    It remains to compare policy improvement. For every \(i,j \in \mathbb N\), the relevant Q-blocks are
    \begin{equation}
        \label{eq:q-control-blocks}
        Q_{i,u_{i}u_{j}}^{(l+1)} =
        \begin{cases}
            R_{ii}+D_{i}^{\top} P_{i}^{(l+1)}D_{i}
            +\displaystyle\int_{\mathcal{E}}F_{i}(e)^{\top}
            P_{i}^{(l+1)}F_{i}(e)\lambda(\mathrm{d}e), & j = i,    \\[2mm]
            D_{i}^{\top} P_{i}^{(l+1)}D_{j}
            +\displaystyle\int_{\mathcal{E}}F_{i}(e)^{\top}
            P_{i}^{(l+1)}F_{j}(e)\lambda(\mathrm{d}e), & j \neq i,
        \end{cases}
    \end{equation}
    and
    \begin{equation}
        \label{eq:q-state-control-blocks}
        Q_{i,u_{i} x}^{(l+1)}
        = B_{i}^{\top} P_{i}^{(l+1)}+D_{i}^{\top} P_{i}^{(l+1)}C
        +\int_{\mathcal{E}}F_{i}(e)^{\top} P_{i}^{(l+1)}E(e)\lambda(\mathrm{d}e).
    \end{equation}
    Comparing~\eqref{eq:q-control-blocks}--\eqref{eq:q-state-control-blocks} with~\eqref{eq:pi-improvement-matrix}--\eqref{eq:pi-improvement-right-hand-side} yields
    \begin{equation*}
        \label{eq:q-improvement-block-equivalence}
        \begin{bmatrix}
            Q_{1,u_{1}u_{1}}^{(l+1)} & \cdots & Q_{1,u_{1}u_{N}}^{(l+1)} \\
            \vdots                   & \ddots & \vdots                   \\
            Q_{N,u_{N}u_{1}}^{(l+1)} & \cdots & Q_{N,u_{N}u_{N}}^{(l+1)}
        \end{bmatrix}
        = \mathcal{M}(P_{1}^{(l+1)},\ldots,P_{N}^{(l+1)}),
        \qquad
        \begin{bmatrix}Q_{1,u_{1}x}^{(l+1)}\\\vdots\\Q_{N,u_{N}x}^{(l+1)}\end{bmatrix}
        = \mathcal{N}(P_{1}^{(l+1)},\ldots,P_{N}^{(l+1)}).
    \end{equation*}
    Substitution into~\eqref{eq:q-gain-update} gives
    \begin{equation*}
        \label{eq:q-gain-update-equivalence}
        \begin{bmatrix}K_{1}^{(l+1)}\\\vdots\\K_{N}^{(l+1)}\end{bmatrix}
        = -\mathcal{M}(P_{1}^{(l+1)},\ldots,P_{N}^{(l+1)})^{-1}
        \mathcal{N}(P_{1}^{(l+1)},\ldots,P_{N}^{(l+1)}),
    \end{equation*}
    which is exactly~\eqref{eq:pi-gain-update}. Therefore, the two algorithms produce the same value matrices and the same updated feedback profile at iteration \(l\).
\end{proof}

\begin{theorem}[Convergence of Algorithm~\ref{alg:off-policy-q-learning}]
    \label{thm:q-convergence}
    Suppose Assumptions~\ref{ass:invariant-ball} and~\ref{ass:quantitative-excitation} hold. If Algorithms~\ref{alg:off-policy-q-learning} and~\ref{alg:policy-iteration} start from the same admissible profile \(K^{(0)}\), then they have the same iterative stability and geometric convergence. The common limit \(K^*\) is the closed-loop Nash equilibrium obtained in Theorem~\ref{thm:pi-convergence}.
\end{theorem}

\begin{proof}
    At \(l = 0\), the two algorithms start from the same feedback profile by assumption, so Lemma~\ref{lem:q-population-equivalence} yields identical \(P_{i}^{(1)}\) and identical updated gains \(K^{(1)}\). Now suppose the feedback profiles still agree at iteration \(l\). Applying the same lemma once more gives equality of every \(P_{i}^{(l+1)}\) as well as of \(K^{(l+1)}\). Thus, by induction, the two algorithms produce exactly the same sequence for all \(l \geq 0\). The stability, geometric convergence, and the closed-loop Nash property therefore follow directly from Theorem~\ref{thm:pi-convergence}.
\end{proof}

\section{Numerical Simulation}
\label{sec:numerical-simulation}
In this section, we provide a numerical example to compare Algorithms~\ref{alg:policy-iteration} and~\ref{alg:off-policy-q-learning}, denoted by A1 and A2, respectively. We take \(N = 2\), \(n = 2\), and \(m_{1} = m_{2} = 1\). The system matrices and the Poisson intensity are
\begin{equation*}
    \label{eq:num-system-coefficients}
    \begin{aligned}
        A & =\begin{bmatrix}-2.09 & -1.78\\-1.77 & -1.89\end{bmatrix}, &
        B = [B_{1}\ B_{2}] & =\begin{bmatrix}0.22 & 0.20\\0.14 & 0.23\end{bmatrix},\\
        C & =\begin{bmatrix}0.67 & 0\\0 & 0.45\end{bmatrix}, &
        D = [D_{1}\ D_{2}] & =\begin{bmatrix}0.08 & 0\\0 & 0.10\end{bmatrix},\\
        E & =\begin{bmatrix}-0.39 & -0.16\\-0.17 & -0.45\end{bmatrix}, &
        F = [F_{1}\ F_{2}] & =\begin{bmatrix}0.11 & 0\\0 & 0.10\end{bmatrix},
        \qquad \lambda = 0.7.
    \end{aligned}
\end{equation*}
The state and control weights are
\begin{equation*}
    \label{eq:num-state-weights}
    \begin{aligned}
        N_{1} & =\begin{bmatrix}6.8 & 2.3\\2.3 & 9.4\end{bmatrix}, &
        N_{2} & =\begin{bmatrix}7.7 & 1.9\\1.9 & 6.8\end{bmatrix},
    \end{aligned}
\end{equation*}
\begin{equation*}
    \label{eq:num-control-weights}
    R_{11} = 0.95,\quad R_{22} = 0.75,\quad R_{12} = 0.45,\quad R_{21} = 0.37.
\end{equation*}
Both algorithms start from \(K_{1}^{(0)}=K_{2}^{(0)}=[0,0]\) with stopping tolerance \(10^{-8}\). This feedback is mean-square stabilizing, but fails the local certificate since \(\|D\mathcal T(0)\|_*\approx4.9928>1\). Nevertheless, the sampled iterations numerically converge. Algorithm~\ref{alg:policy-iteration} terminates after \(15\) updates, yielding the following final value matrices and feedback gains
\begin{equation*}
    \label{eq:num-policy-iteration-results}
    \begin{aligned}
        P_{1,\mathrm{A1}} & =\begin{bmatrix}
            12.3359 & -11.6767 \\-11.6767 & 14.0121
        \end{bmatrix}, &
        K_{1,\mathrm{A1}} & =\begin{bmatrix}-1.2650 & 0.5337\end{bmatrix},\\
        P_{2,\mathrm{A1}} & =\begin{bmatrix}
            10.8516 & -9.9738 \\-9.9738 & 11.6903
        \end{bmatrix}, &
        K_{2,\mathrm{A1}} & =\begin{bmatrix}0.4856 & -0.9275\end{bmatrix}.
    \end{aligned}
\end{equation*}
Algorithm~\ref{alg:off-policy-q-learning} terminates after \(16\) updates and yields
\begin{equation*}
    \label{eq:num-q-learning-results}
    \begin{aligned}
        P_{1,\mathrm{A2}} & =\begin{bmatrix}
            12.2994 & -11.6176 \\-11.6176 & 13.9214
        \end{bmatrix}, &
        K_{1,\mathrm{A2}} & =\begin{bmatrix}-1.2813 & 0.5551\end{bmatrix},\\
        P_{2,\mathrm{A2}} & =\begin{bmatrix}
            10.7511 & -9.8508 \\-9.8508 & 11.5368
        \end{bmatrix}, &
        K_{2,\mathrm{A2}} & =\begin{bmatrix}0.4786 & -0.9159\end{bmatrix}.
    \end{aligned}
\end{equation*}
We then go back to check the results of the entire feedback iteration sequence. Every saved iterate is mean-square stabilizing. Starting from the third update, all iterates stay inside closed balls centered at \(K_{\mathrm{A1}}^{(3)}\) and \(K_{\mathrm{A2}}^{(3)}\), with radii \(0.06\) and \(0.09\), respectively. We check that Assumption~\ref{ass:local-certificate} holds on both balls using interval bounds, so Assumption~\ref{ass:invariant-ball} is satisfied. As a result, Theorem~\ref{thm:pi-convergence} ensures convergence if we run exact policy iteration starting from either of the third-update gains. The sampled numerical sequences support this local convergence behavior.

\begin{figure*}[t]
	\centering
	\begin{minipage}[t]{0.32\textwidth}
		\centering
		\includegraphics[width=\linewidth]{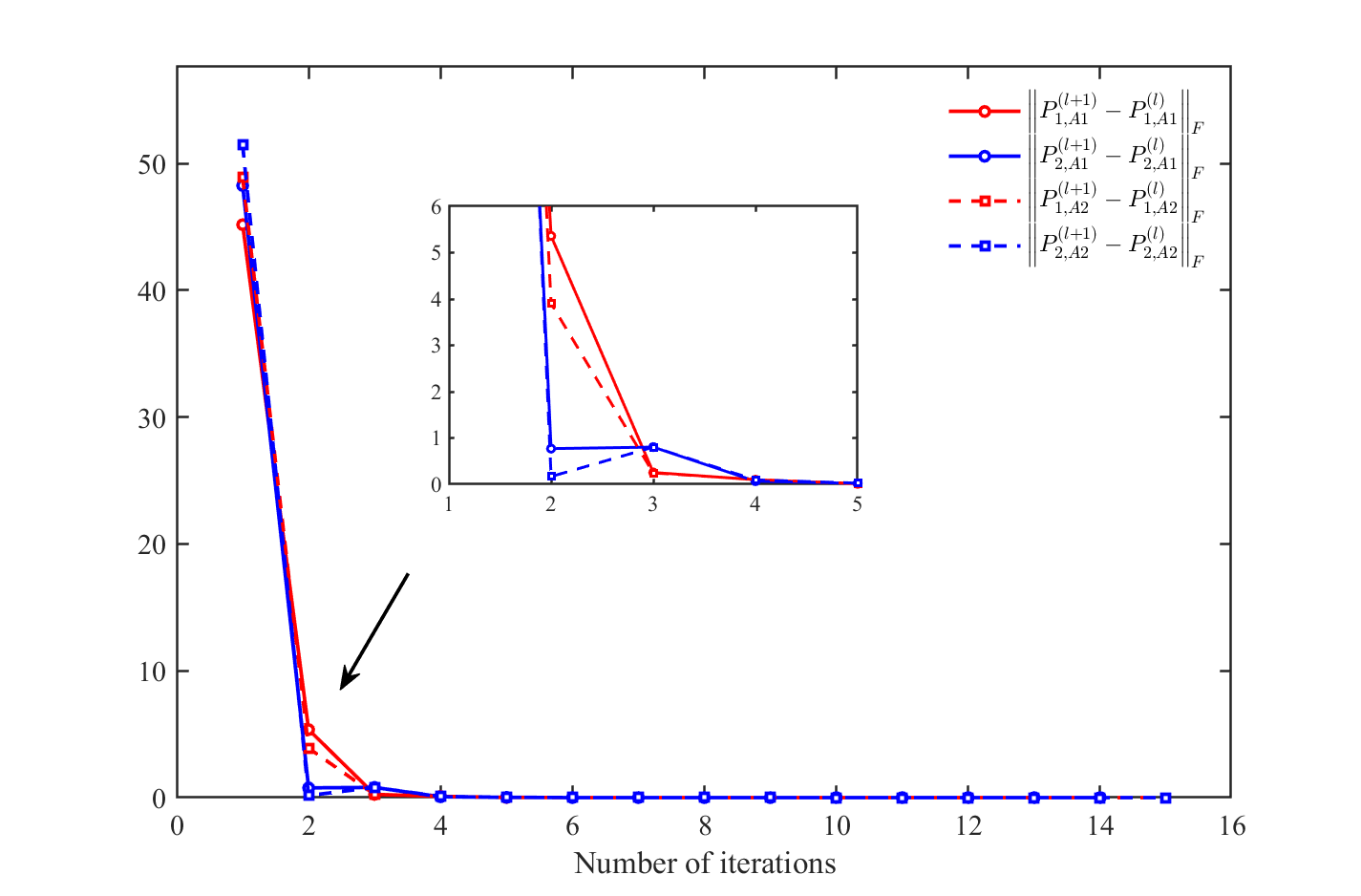}\\
		\small (a) 
	\end{minipage}
	\hfill
	\begin{minipage}[t]{0.32\textwidth}
		\centering
		\includegraphics[width=\linewidth]{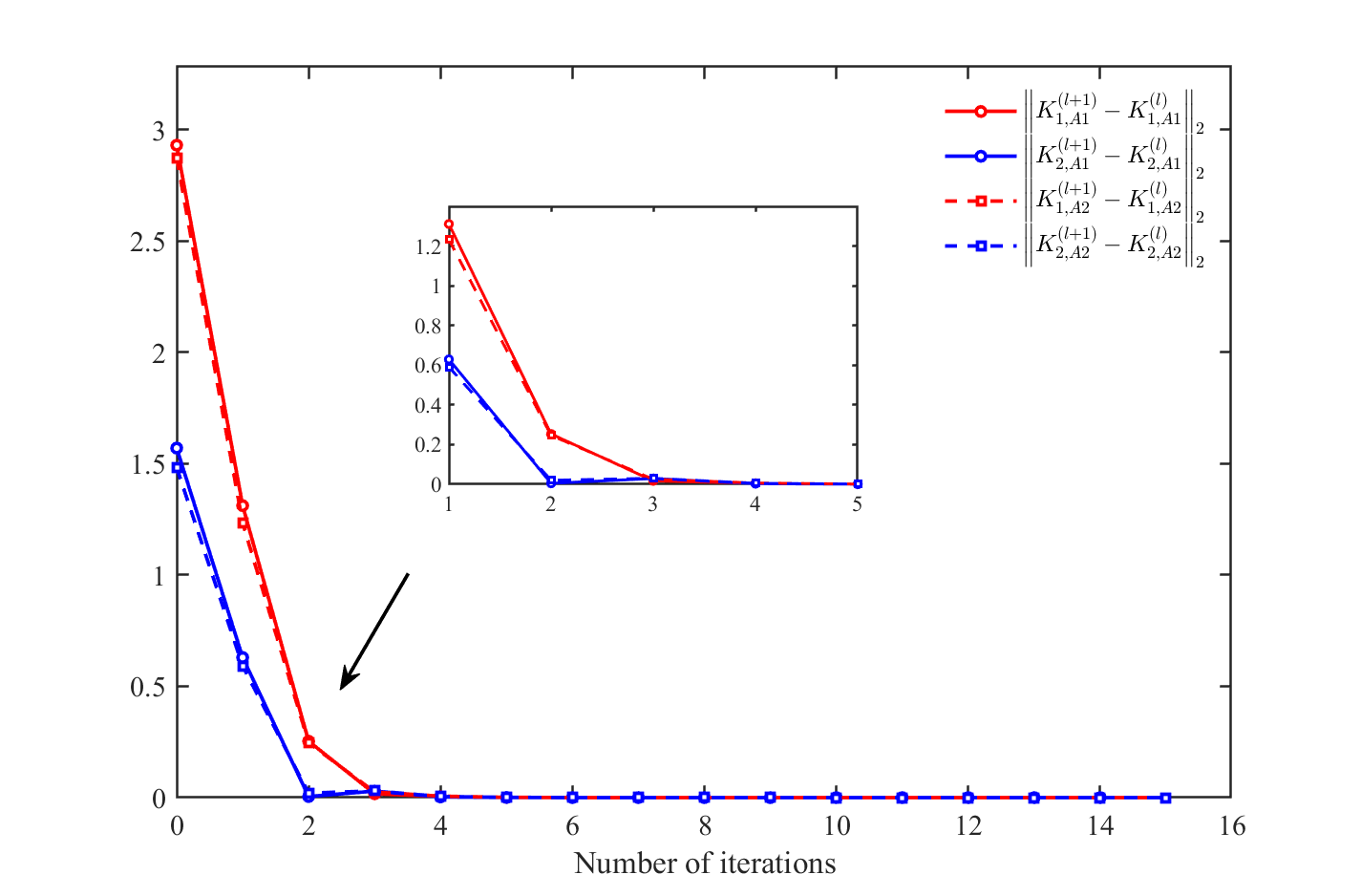}\\
		\small (b) 
	\end{minipage}
	\hfill
	\begin{minipage}[t]{0.32\textwidth}
		\centering
		\includegraphics[width=\linewidth]{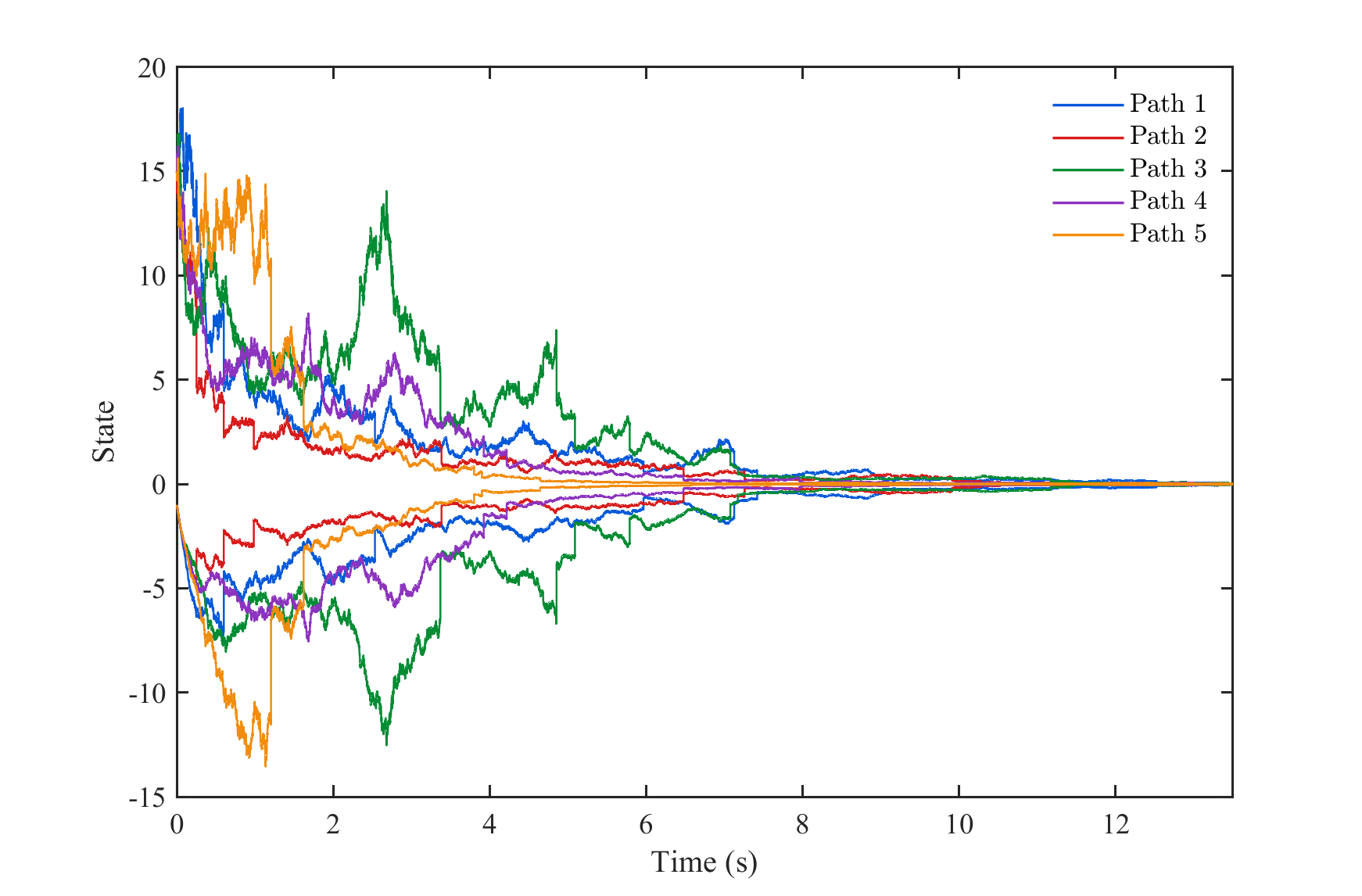}\\
		\small (c)
	\end{minipage}
	\caption{(a) Convergence of \(P_i^{(l)}\) for Algorithms~\ref{alg:policy-iteration} and~\ref{alg:off-policy-q-learning};
		(b) Convergence of \(K_i^{(l)}\) for Algorithms~\ref{alg:policy-iteration} and~\ref{alg:off-policy-q-learning};
		(c) Closed-loop state sample paths under the final feedback profile.}
	\label{fig:simulation-results}
\end{figure*}

To compare the final results of Algorithms~\ref{alg:policy-iteration} and~\ref{alg:off-policy-q-learning}, we calculate the Frobenius-norm differences between their final iterates. The individual differences are
\begin{equation*}
    \label{eq:num-individual-differences}
    \begin{aligned}
        \|P_{1,\mathrm{A1}}-P_{1,\mathrm{A2}}\|_F & =0.1286, &
        \|P_{2,\mathrm{A1}}-P_{2,\mathrm{A2}}\|_F & =0.2529,\\
        \|K_{1,\mathrm{A1}}-K_{1,\mathrm{A2}}\|_F & =0.0270, &
        \|K_{2,\mathrm{A1}}-K_{2,\mathrm{A2}}\|_F & =0.0136.
    \end{aligned}
\end{equation*}
Consequently, the total differences between the two algorithms are
\begin{equation*}
    \label{eq:num-total-differences}
    \begin{aligned}
        \sum_{i = 1}^{2} \|P_{i,\mathrm{A1}}-P_{i,\mathrm{A2}}\|_F & =0.3815,\\
        \sum_{i = 1}^{2} \|K_{i,\mathrm{A1}}-K_{i,\mathrm{A2}}\|_F & =0.0405.
    \end{aligned}
\end{equation*}
The differences between the two algorithms are therefore very small. Fig.~\ref{fig:simulation-results}(a) and (b) show the Frobenius norms of the value-matrix increments and the Euclidean norms of the feedback-gain increments, respectively, for both algorithms. Fig.~\ref{fig:simulation-results}(c) shows five independent state sample paths under the final feedback profile, starting from \(x_0=(15,-1)^\top\).

Numerical results show that both algorithms converge rapidly, maintain the mean-square stability of the system during iterations, and yield close limiting matrices and feedback gains. In particular, Algorithm~\ref{alg:off-policy-q-learning} reuses the same behavior-policy data across iterations without requiring the system coefficients or jump intensity.

\section{Conclusion}
\label{sec:conclusion}
In this paper, we developed two medel-free reinforcement learning methods for infinite-horizon \(N\)-player linear-quadratic nonzero-sum stochastic differential games with controlled jump diffusions. The methods consider the coupling between players' feedback gains induced by control-dependent diffusion and jump amplitudes. A verifiable certificate on a closed ball centered at an admissible initial feedback profile guarantees mean-square stability of all policy iterates and local geometric convergence to a closed-loop Nash equilibrium. Totally model-free off-policy Q-learning reuses the same behavior-policy data across iterations. Under quantitative persistence of excitation, exact population regression reproduces the policy-iteration sequence from the same initial profile. In the two-player example, although the zero initial feedback fails the certificate, exact policy iteration is guaranteed to converge when started from either algorithm’s third iterate. Future research may extend the present framework to nonlinear systems and large-scale networked games.

\bibliographystyle{plain}
\bibliography{ref}

\end{document}